\documentclass[reqno, 10pt]{amsart}

\usepackage[papersize={18cm, 26cm}, margin=2.65cm, tmargin=1.8cm, bmargin=1.4cm]{geometry}
\usepackage{amsmath}
\usepackage{amssymb}
\usepackage{amsmath}
\usepackage{amssymb}
\usepackage{amssymb}
\usepackage{mathrsfs}
\usepackage{amsmath}
\usepackage{amsmath,amsthm,amssymb,amscd}
\usepackage{latexsym}
\usepackage{color}
\usepackage{graphicx}
\usepackage[colorlinks=true]{hyperref}
\hypersetup{urlcolor=blue, citecolor=blue,linkcolor=blue}

\newtheorem{theorem}{Theorem}[section]
\newtheorem{corollary}[theorem]{Corollary}

\newtheorem{remark}[theorem]{Remark}
\newtheorem{example}[theorem]{Example}

\theoremstyle{definition} \theoremstyle{remark}
\def\C{\mathbb C}
\def\R{\mathbb R}

\def\Z{\mathbb Z}

\def\Z{\mathbb Z}

\def\<{\langle}
\def\>{\rangle}

\def\diag {\mathrm{diag}}

\newcommand{\tr}{\operatorname{tr}}

\def\F{\mathbb{F}}
\renewcommand{\i}{\mathbf{i}}
\newcommand{\be}[2]{\begin{#1} #2 \end{#1}}  %% environment
\newcommand{\mtx}[1]{\begin{pmatrix} #1 \end{pmatrix}} %% matrix
\newcommand{\wt}[1]{\widetilde{#1}}

\newcommand{\ol}[1]{\overline{#1}}

\begin{document}

\title{Linear k-power preservers and trace of power-product preservers}

\author[H. Huang] {Huajun Huang}
\address{Department of Mathematics and Statistics, Auburn University,
AL 36849--5310, USA}
  \email{huanghu@auburn.edu}

\author[M. Tsai] {Ming-Cheng Tsai}
\address{General Education Center,
Taipei University of Technology, Taipei 10608, Taiwan}
\email{mctsai2@mail.ntut.edu.tw}

\begin{abstract}
Let $V$ be  the set  of $n\times n$ complex or real general matrices, Hermitian matrices, symmetric matrices, positive definite (resp. semi-definite) matrices,  diagonal matrices, or upper triangular matrices.  
Fix    $k\in\Z\setminus \{0, 1\}$.  We characterize linear maps  $\psi:V\to V$  that satisfy $\psi(A^k)=\psi(A)^k$  on an open neighborhood $S$ of $I_n$ in $V$.
The $k$-power preservers are necessarily $k$-potent preservers, and $k=2$ corresponds to Jordan homomorphisms. 
 Applying the   results, we characterize maps $\phi,\psi:V\to V$ that satisfy
``$ \tr(\phi(A)\psi(B)^k)=\tr(AB^k)$ for all $A\in V$, $B\in S$, and $\psi$ is linear'' or ``$ \tr(\phi(A)\psi(B)^k)=\tr(AB^k)$ for all $A, B\in S$ and both $\phi$ and $\psi$ are  linear.''  
The characterizations systematically extend existing results in literature, and they have many applications in areas like quantum information theory. 
 Some structural theorems and power series over matrices  are widely used in our characterizations. 
\end{abstract}

%These results generalize the corresponding $k$-power preserving and trace preserving problems. 

%\keywords{trace of products, preserver, ?}

%\footer{Mathematics Subject Classification 2020: Primary 15A21, Secondary 15A23, 14L30.}

\maketitle

\section{Introduction}

Preserver problem is one of the most active research areas in matrix theory (e.g. \cite{Li01, Li92, Molnar, Pierce}). Researchers would like to characterize the maps on a given space of matrices preserving certain subsets, functions or relations. One of the preserver problems concerns maps $\psi$ on some sets $V$ of matrices which preserves $k$-power for a fixed integer $k\geq 2$, that is, $\psi(A^k)=\psi(A)^k$ for any $A\in V$ (e.g. \cite{ChanLim, Molnar, ZTC}). 
The  $k$-power preservers form a special class of polynomial preservers.   
One important reason of this problem lies on the fact that the case $k=2$ corresponds to  Jordan homomorphisms. 
Moreover, every $k$-power preserver is also a $k$-potent preserver,
that is, $A^k=A$ imply that $\psi(A)^k=\psi(A)$ for any $A\in V$. Some researches on $k$-potent preservers can be found in \cite{BS93, ZC06, ZTC}.

Given a field $\F$, let $\mathcal{M}_n(\F)$, $\mathcal{S}_n(\F)$, $\mathcal{D}_n(\F)$,   $\mathcal{N}_n(\F)$, and $\mathcal{T}_n(\F)$ denote the set of $n\times n$ general, symmetric, diagonal, strictly upper triangular, and  upper triangular matrices over $\F$, respectively. When $\F$ is the complex field $\C$, we may write $\mathcal{M}_n$ instead of $\mathcal{M}_n(\C)$, and so on. Let $\mathcal{H}_n$, $\mathcal{P}_n$, and $\ol{\mathcal{P}_n}$ denote the set of complex Hermitian, positive definite, and positive semidefinite matrices, and $\mathcal{H}_n(\R)=\mathcal{S}_n(\R)$, $\mathcal{P}_n(\R)$, and $\ol{\mathcal{P}_n(\R)}$ the corresponding set of real matrices, respectively.
 A matrix space is a subspace of  $\mathcal{M}_{m,n}(\F)$  for certain $m,n\in\Z_+$.   Let $A^t$ (resp. $A^*$) denote the transpose (resp. conjugate transpose) of a matrix $A$.

In 1951, Kadison (\cite{Kadison}) showed that a Jordan $*$-isomorphism on $\mathcal{M}_n$, namely, a bijective linear map with $\psi(A^2)=\psi(A)^2$ and $\psi(A^*)=\psi(A)^*$ for all $A\in \mathcal{M}_n$, is the direct sum of a $*$-isomorphism and a $*$-anti-isomorphism. Hence $\psi(A)=UAU^*$ for all $A\in \mathcal{M}_n$ or $\psi(A)=UA^{T}U^*$ for all $A\in \mathcal{M}_n$ by \cite[Theorem A.8]{Molnar}. Let  $k\geq2$ be  a fixed integer. In 1992, Chan and Lim (\cite{ChanLim}) determined a nonzero linear operator $\psi: \mathcal{M}_n(\F)\to \mathcal{M}_n(\F)$ (resp. $\psi: \mathcal{S}_n(\F)\to \mathcal{S}_n(\F)$)  such that $\psi(A^k)=\psi(A)^k$ for all $A\in \mathcal{M}_n(\F)$ (resp. $\mathcal{S}_n(\F)$) (See Theorems \ref{thm: Mn ChanLim} and  \ref{thm: Sn ChanLim}). In 1998, Bre\u{s}ar, Martindale, and Miers considered additive maps of general prime rings to solve an analogous problem by using the deep algebraic techniques (\cite{BMM}). Monl\'{a}r \cite[P6]{Molnar} described a particular case of their result which extends Theorem \ref{thm: Mn ChanLim} to surjective linear operators on $\mathcal{B(H)}$. 
In 2004, Cao and Zhang determined additive $k$-power preserver on $\mathcal{M}_n(\F)$ and $\mathcal{S}_n(\F)$ (\cite{Cao04}). They also  characterized injective additive $k$-power preserver on $\mathcal{T}_n(\F)$ (\cite{Cao05} or \cite[Theorem 6.5.2]{ZTC}), which leads to injective linear $k$-power preserver on $\mathcal{T}_n(\F)$ (see Theorem \ref{thm: Tn injective k power preserver}).
In 2006, Cao and Zhang also characterized linear $k$-power preservers from $\mathcal{M}_n(\F)$ to $\mathcal{M}_m(\F)$ and from $\mathcal{S}_n(\F)$ to $\mathcal{M}_m(\F)$ (resp. $\mathcal{S}_m(\F)$) \cite{ZC06}.

In this article, given an integer $k \in\Z\setminus\{0,1\}$, we show that a unital linear map $\psi:V\to W$ between matrix spaces preserving $k$-powers on a neighborhood of identity  must preserve all integer powers (Theorem \ref{thm: k-power all powers}). Then we characterize,
for $\F=\C$ and $\R$, linear operators  on  sets $V=
\mathcal{M}_n(\F)$, $\mathcal{H}_n$, $\mathcal{S}_n(\F)$,  $\mathcal{P}_n$, $\mathcal{P}_n(\R)$, and $\mathcal{D}_n (\F)$ that satisfy $\psi(A^k)=\psi(A)^k$  on an open neighborhood of $I_n$ in $V$. In the following descriptions,  $P\in \mathcal{M}_n(\F)$ is invertible, $U\in\mathcal{M}_n(\F)$ is unitary, 
$O\in\mathcal{M}_n(\F)$ is orthogonal,  and $\lambda\in\F$ satisfies that $\lambda ^{k-1}=1$.
\begin{enumerate}
\item
$V=\mathcal{M}_n(\F)$ (Theorem \ref{thm: Mn k-power negative}):
$\psi(A)=\lambda PAP^{-1}$ or $\psi(A)=\lambda PA^t P^{-1}$. 

\item
$V=\mathcal{H}_n$ (Theorem \ref{thm: Hn linear map preserving k power}): 
When $k$ is even, $\psi(A)=U^*AU$ or $\psi(A)=U^*A^t U$.
When $k$ is odd, $\psi(A)=\pm U^*AU$ or $\psi(A)=\pm U^*A^t U$.

\item
$V=\mathcal{S}_n(\F)$ (Theorem \ref{thm: Sn k-power negative}): 
$\psi(A)=\lambda OAO^{t}$. 

\item
$V=\mathcal{P}_n$ or $\mathcal{P}_n(\R)$ (Theorem \ref{thm: Pn linear map preserving k power}):
$\psi(A)=U^*AU$ or $\psi(A)=U^*A^t U$.

\item
$V=\mathcal{D}_n(\F)$ (Theorem \ref {thm: Dn k power}):
$\psi(A)=\psi(I_n)\diag(  f_{p(1)}(A),\ldots, f_{p(n)}(A))$, in which
$\psi(I_n)^k=\psi(I_n)$, $p:\{1,\ldots,n\}\to\{0,1,\ldots,n\}$ is a function, and $f_i:\mathcal{D}_n(\F)\to\F$  $(i=0,1,\ldots,n)$ satisfy that,
for $A=\diag(a_1,\ldots,a_n)$,
$f_0(A)=0$ and $f_i(A)=a_i$ for $i=1,\ldots,n$. 

\item
$V=\mathcal{T}_n(\F)$ (Theorem \ref{thm: Tn injective k-power negative} for $n\ge 3$): 
$\psi(A)=\lambda PAP^{-1}$ or $\psi(A)=\lambda PA^{-}P^{-1}$, in which $P\in \mathcal{T}_n(\F)$
and $A^{-}=(a_{n+1-j,n+1-i})$ if $A=(a_{ij})$. 

\end{enumerate}
Our results on $\mathcal{M}_n(\F)$ and $\mathcal{S}_n(\F)$ extend Chan and Lim's results in Theorems \ref{thm: Mn ChanLim} and  \ref{thm: Sn ChanLim},
and  result on $\mathcal{T}_n(\F)$ extend  Cao and Zhang's linear version result in \cite{Cao05}.

Another topic is the study of a linear map $\phi$ from a matrix set $S$ to another matrix set $T$ preserving trace equation. In 1931, Wigner's unitary-antiunitary theorem {\cite[p.12]{Molnar}} says that if $\phi$ is a bijective map  defined on the set of all rank one projections on a Hilbert space $H$ satisfying
%\begin{equation*}%\label{wigner}
\begin{equation}\label{tr-prod-2-preserving}
\tr(\phi(A)\phi(B))=\tr(AB),
\end{equation}
%\end{equation*}
%that is, $\phi$ preserves the transition probability between the pure states $P$ and $Q$ (or the angle between the ranges of $P$ and $Q$),
then there is an either unitary or antiunitary operator $U$ on $H$ such that $\phi(P)=U^{*}PU$ or $\phi(P)=U^{*}P^{t}U$ for all rank one projections $P$. In 1963, Uhlhorn generalized Wigner's theorem to show that the same conclusion holds if the equality $\tr(\phi(P)\phi(Q))=\tr(PQ)$ is replaced by $\tr(\phi(P)\phi(Q))=0 \Leftrightarrow \tr(PQ)=0$ (see \cite{Uh63}).
% The study of this topic has been applied to areas like   invertible quantum channels.

%Let $\mathcal{M}_{m,n}$ denote the set of $m\times n$ complex matrices.

In 2002, Moln\'{a}r (in the proof of \cite[Theorem 1]{Molnar2}) showed that maps $\phi$ on the space of all bounded linear operators on a Banach space $B(X)$ satisfying \eqref{tr-prod-2-preserving} for $A\in B(X)$, rank one operator $B\in B(X)$ are linear. In 2012, Li, Plevnik, and {\v S}emrl \cite{Li12} characterized bijective maps $\phi: S\rightarrow S$ satisfying $\tr(\phi(A)\phi(B))=c  \Leftrightarrow \tr(AB)=c$ for a given real number $c$, where $S$ is $\mathcal{H}_n$, $\mathcal{S}_n(\R)$,  or the set of  rank one projections.

In \cite[Lemma 3.6]{Huang16}, Huang et al showed that the following statements are equivalent for
a unital map $\phi$ on $\mathcal{P}_n$:
\be{enumerate}{
\item
$\tr(\phi(A)\phi(B))=\tr(AB)$ for $A,B\in \mathcal{P}_n$;
\item
$\tr(\phi(A)\phi(B)^{-1})=\tr(AB^{-1})$ for $A, B\in \mathcal{P}_n$;
\item
$\phi(A)=U^*AU$ or $U^*A^t U$ for a unitary matrix $U$.
}
The authors  also determined the cases if $\phi$ is not assuming unital, the set $\mathcal{P}_n$ is replaced by another set like $\mathcal{M}_n$,  $\mathcal{S}_n$ , $\mathcal{T}_n$, or $\mathcal{D}_n$. In \cite[Theorem 3.8]{Leung16}, Leung, Ng, and Wong considered the relation \eqref{tr-prod-2-preserving} on infinite dimensional space.

%\iffalse
Let $\langle S\rangle$ denote the subspace spanned by a subset $S$ of a vector space. 
Recently,   Huang and Tsai studied two maps preserving  trace of product \cite{Huang21}. Suppose two maps $\phi: V_1\to W_1$ and $\psi: V_2\to W_2$ between subsets of matrix spaces over a field $\F$ under some conditions satisfy
 \begin{equation}\label{tr-prod-2d-preserving}
\tr(\phi(A)\psi(B))=\tr(AB)
\end{equation}
for all $A\in V_1$, $B\in V_2$. The authors showed that  these two maps can be extended to bijective linear maps $\widetilde\phi : \langle V_1\rangle \to  \langle W_1\rangle$ and $\widetilde\psi: \langle V_2\rangle \to  \langle W_2\rangle$ that satisfy
 $\tr(\widetilde\phi(A)\widetilde\psi(B))=\tr(AB)$
for all $A\in \langle V_1\rangle$, $B\in \langle V_2\rangle$ (see Theorem \ref{thm: two maps preserving trace}).  Hence  when a matrix space $V$ is closed under conjugate transpose,  every  linear bijection $\phi:V\to V$  corresponds to a unique
 linear bijection $\psi:V\to V$ that makes \eqref{tr-prod-2d-preserving} hold (see Corollary \ref{thm: trace preserver on special spaces}). %\cite[Corollary 2.2]{Huang21}. 
 Therefore, each of $\phi$ and $\psi$ has no specific form. 
 
One natural question is to ask when the following equality holds for a fixed $k\in\Z\setminus\{0, 1\}$: 
 \begin{equation}\label{tr-prod-k-preserving}
 \tr(\phi(A)\psi(B)^k)=\tr(AB^k).
\end{equation}
 In this article, we   use our characterization of linear  $k$-power preservers on an open neighborhood $S$ of $I_n$ in   $V$ to  characterize maps $\phi,\psi:V\to V$   under one of the assumptions: 
\begin{enumerate}
\item equality \eqref{tr-prod-k-preserving} holds for all $A\in V$, $B\in S$, and $\psi$ is linear, or
\item equality \eqref{tr-prod-k-preserving} holds for all $A, B\in S$ and both $\phi$ and $\psi$ are  linear,
\end{enumerate}
for the sets $V=\mathcal{M}_n$, $\mathcal{H}_n$, $\mathcal{P}_n$, $\mathcal{S}_n$, $\mathcal{D}_n$, and their real counterparts. In the following descriptions,   $\F=\C$ or $\R$,  $P, Q\in M_n(\F)$ are invertible, $U\in M_n(\F)$ is unitary, $O\in M_n(\F)$ is orthogonal, and $c\in\F\setminus\{0\}$.

%The map $\phi : \mathcal{S}\to \mathcal{S}$ and the linear map $\psi : \mathcal{S}\to \mathcal{S}$  that satisfies \eqref{1k} are described as follow:
\begin{enumerate}

\item   $V=\mathcal{M}_n(\F)$ (Theorem \ref{Two traces}): 
\begin{enumerate}
\item When $k=-1$, $\phi(A)=PAQ$ and $\psi(B)=PBQ$, or $\phi(A)=PA^t Q$ and $\psi(B)=PB^t Q$.
\item 
When $k\in\Z\setminus\{-1,0,1\}$, $\phi(A) = c^{-k} PAP^{-1}$ and 
$\psi(B) = cPBP^{-1}$, or $\phi(A) = c^{-k} PA^t P^{-1}$ and
$\psi(B) = cPB^t P^{-1}$.  
\end{enumerate}

\item   $V=\mathcal{H}_n$ (Theorem \ref{thm: Hn Two traces}): 
\begin{enumerate}
\item When $k=-1$, 
$\phi(A) = cP^*AP$ and
$\psi(B) = cP^*BP$, or
$\phi(A) = cP^*A^t P$ and
$\psi(B) = cP^*B^t P$,
for $c\in\{1,-1\}$. 
\item When $k\in \Z\setminus\{-1,0,1\}$,
$\phi(A) = c^{-k} U^*AU$ and
$\psi(B) = cU^*BU$, or 
$\phi(A) = c^{-k} U^*A^t U$ and
$\psi(B) = cU^*B^t U$, for $c\in\R\setminus\{0\}$. 
 \end{enumerate}

%\mathcal{S}

\item    $V=\mathcal{S}_n(\F)$  (Theorem \ref{thm: Sn Two traces}): 
\begin{enumerate}
\item 
When $k=-1$, $\phi(A)=cP A P^t$ and $\psi(B)= cP B P^t$.
\item
When $k\in\Z\setminus\{-1,0,1\}$, $\phi(A) = c^{-k} O AO^t$ and $\psi(B) = cOBO^t$. 
\end{enumerate}
 
\item $V=\mathcal{P}_n$ and $\mathcal{P}_n(\R)$  (Theorem \ref{thm: Pn preserve trace power open set}): 
$\phi(A) = c^{-k} U^*AU$ and 
$\psi(B) = cU^*BU$, 
or 
$\phi(A) = c^{-k} U^*A^t U$ and
$\psi(B) = cU^*B^t U$, in which $c\in\R^+$. 
Characterizations under some other assumptions are also given as special cases
of Theorem \ref{thm: Pn preserve trace power} (Huang, Tsai \cite{Huang21}). 
%% Theorem \ref{thm: Pn preserve trace power} (\cite{Huang21}) shows that given $a, b, c, d\in\R\setminus\{0\}$,  two maps    $\phi, \psi:\mathcal{P}_n\to \ol{\mathcal{P}_n}$   satisfy  $\tr (\phi(A)^{a} \psi(B)^{b})=\tr(A^c B^d)$ for $A, B\in \mathcal{P}_n$  if and only if $\phi(A)=(P^*A^c P)^{1/a}$ and $\psi(B)=(P^{-1} B^d P^{-*})^{1/b}$, or  $\phi(A)=[P^*(A^t)^{c}P]^{1/a}$ and $\psi(B)=[P^{-1} (B^t)^{d} P^{-*}]^{1/b}$.    All results have counterparts for maps  $\phi, \psi:\mathcal{P}_n(\R)\to \ol{\mathcal{P}_n(\R)}$. 

\item    $V=\mathcal{D}_n(\F)$  (Theorem \ref{thm: Dn Two traces2}): $\phi (A)=PC^{-k}A P^{-1},\ \psi (B)=PC B P^{-1}$ where $P$ is a permutation matrix  and $C=\mathcal{D}_n(\F)$ is diagonal and invertible.

\item $V=\mathcal{T}_n(\F)$  (Theorem \ref{Tn Two traces}): 
$\phi$ and $\psi$  send $\mathcal{N}_n(\F)$ 
to $\mathcal{N}_n(\F)$, $(\operatorname{D}\circ\phi)|_{\mathcal{D}_n(\F)}$ and $(\operatorname{D}\circ\psi)|_{\mathcal{D}_n(\F)}$
are characterized by Theorem \ref{thm: Dn Two traces2},  and
$\operatorname{D}\circ\phi=\operatorname{D}\circ\phi\circ \operatorname{D}$. Here 
 $\operatorname{D}$ denotes the map that sends $A\in \mathcal{T}_n(\F)$ to the diagonal matrix with the same diagonal as $A$.

\end{enumerate}

The sets  $\mathcal{M}_n$, $\mathcal{H}_n$, $\mathcal{P}_n$, $\mathcal{S}_n,$ $\mathcal{D}_n$, and their real counterparts are closed under conjugate transpose.  
In these sets, $\tr(AB)=\langle A^*, B\rangle$ for the standard inner product. 
Our trace of product preservers   can also be interpreted as inner product preservers, which have wide applications in research areas like quantum information theory.

\section{Preliminary}

\subsection{Linear operators  preserving powers}

We show below that: given $k\in\Z\setminus\{0,1\}$, a unital  linear map $\psi: V\to W$ between matrix spaces preserving $k$-powers on a neighborhood of identity in $V$ must preserve all   integer powers.
Let $\Z_+$ (resp. $\Z_-$) denote the set of all positive (resp. negative) integers.

\be{theorem}
{\label{thm: k-power all powers} 
Let $\F=\C$ or $\R$. 
Let $V\subseteq \mathcal{M}_p(\F)$ and $W\subseteq \mathcal{M}_q(\F)$ be   matrix spaces. Fix $k \in\Z\setminus\{0,1\}$.
\begin{enumerate}
\item
Suppose the identity matrix $I_p\in V$  and
  $A^k\in V$ for all matrices $A$ in an open neighborhood $\mathcal{S}_V$ of $I_p$ in $V$ consisting of invertible matrices.
Then
\be{eqnarray}
{
\{AB+BA: A, B\in V\} &\subseteq&  V, \\
\{A^{-1}: A\in V\text{ is invertible}\} &\subseteq&  V.
}
In particular,
\be{eqnarray}{
 \{A^r:A\in V\}\subseteq V, &&  r\in\Z_+, \quad \text{and}
 \\
 \{A^r:A\in V\text{ is invertible}\}\subseteq V, &&  r\in\Z_-.
}

\item
 Suppose $I_p\in V$, $I_q\in W$,  and  $A^k\in V$ for all matrices $A$    in an open neighborhood $\mathcal{S}_V$ of $I_p$ in $V$
 consisting of invertible matrices.
Suppose $\psi:V\to W$ is a linear map that satisfies the following conditions:
\be{eqnarray}
{\label{k-power I}
\psi(I_p) &=& I_q,
\\\label{k-power A^k}
\psi(A^k)&=& \psi(A)^k,\quad   A\in \mathcal{S}_V.
}
Then
\be{eqnarray}{\label{psi(AB+BA) identity}
\psi(AB +BA)&=&\psi(A)\psi(B)+\psi(B)\psi(A),\quad  A,B\in V,
\\\label{psi(A) inverse identity}
\psi(A^{-1}) &=& \psi(A)^{-1},\quad  \text{ invertible } A\in V.
}
In particular,
\be{eqnarray}
{\label{psi(A positive power) identity}
\psi(A^r)=\psi(A)^r,&&    A\in V,\   r\in\Z_+,\quad\text{and}
\\\label{psi(A negative power) identity}
\psi(A^r)=\psi(A)^r,&&    \text{ invertible } A\in V,\   r\in\Z_-.
}
\end{enumerate}
}

\begin{proof}
We prove the complex case. The real case is done similarly.
\begin{enumerate}
\item
For each $A\in V\setminus\{0\}$, there is $\epsilon>0$ such that $I_p+xA\in \mathcal{S}_V$ for all $x\in\C$ with $|x|<\min\{\epsilon, \frac{1}{\|A\|}\}$. Thus
\be{equation}
{\label{I+xA k power}
(I_p+xA)^k=I_p+xkA+x^2 \frac{k(k-1)}{2} A^2+\cdots \in V.
}
The second derivative
\be{equation}{
 \left. \frac{d^2}{dx^2}(I_p+xA)^k\right|_{x=0} =k(k-1)A^2\in V.
}
Since $k\not\in\{0,1\}$, we have $A^2\in V$ for all $A\in V$. Therefore, for $A, B\in V$,
\be{equation}{
AB+BA= (A+B)^2-A^2-B^2\in V.
}
In particular, $A\in V$ implies that $A^r\in V$ for all $r\in\Z_+$.

Cayley-Hamilton theorem implies that every invertible matrix $A$ satisfies that $A^{-1}=f(A)$ for a certain polynomial $f(x)\in\F[x]$.
Therefore, $A^{-1}\in V$, so that $A^{r}\in V$ for all $r\in\Z_-$.

\item Now suppose \eqref{k-power I} and \eqref{k-power A^k} hold.
The proof is proceeded similarly to the proof of part (1).
For every $A\in  V$, there is $\epsilon>0$ such that for all $x\in \C$ with $|x|<\min\{\epsilon,\frac{1}{\|A\|},\frac{1}{\|\psi(A)\|}\}$,
\be{eqnarray}
{ \label{k-power I+xA psi}
(\psi(I_p+xA))^k &=&
 I_q+xk\psi(A)+x^2\frac{k(k-1)}{2}\psi(A)^2+\cdots\in W, \qquad
\\ \label{k-power I+xA}
\psi((I_p+xA)^k)
&=&
I_q+xk\psi(A)+x^2\frac{k(k-1)}{2} \psi(A^2)+\cdots\in W.
}
Since \eqref{k-power I+xA psi} and \eqref{k-power I+xA} equal, we have
\be{equation}
{\label{Monlar A square equality}
\psi(A)^2=\psi(A^2),\quad   A\in  V.
}
Therefore, for $A, B\in V$,
\be{equation}{
\psi((A+B)^2)=\psi(A+B)^2
}
We get \eqref{psi(AB+BA) identity}: $\psi(AB+BA)=\psi(A)\psi(B)+\psi(B)\psi(A)$. In particular  $\psi(A^r)=\psi(A)^r$ for all $A\in V$ and $r\in\Z_+$.

Every invertible $A\in V$ can be expressed as $A^{-1}=f(A)$ for a certain polynomial  $f(x)\in\F[x]$.
Then $\psi(A^{-1})=\psi(f(A))=f(\psi(A))$ is commuting with $\psi(A)$.
Hence
\be{equation}{
2\psi(A^{-1})\psi(A)=\psi(A^{-1})\psi(A)+\psi(A)\psi(A^{-1})=\psi(A^{-1}A+AA^{-1})=2I_q.
}
We get $\psi(A^{-1})=\psi(A)^{-1}$. Therefore, $\psi(A^r)=\psi(A)^r$ for all $r\in\Z_-$.
\qedhere
\end{enumerate}
\end{proof}

%%%%%%%%% DUPLICATED ABOVE %%%%%%%%%%%

Theorem \ref{thm: k-power all powers} is powerful in exploring $k$-power preservers in  matrix spaces.  Note that
every $k$-power preserver is  a $k$-potent preserver. Theorem \ref{thm: k-power all powers}  can also be used to
investigate $k$-potent preservers in  matrix spaces.

\subsection{Two maps preserving trace of product}

We recall two  results about two maps preserving   trace of product in \cite{Huang21}.  They are handy in proving linear bijectivity
of  maps preserving   trace of products. Recall that if $S$ is a subset of a vector space, then
$\langle S\rangle$ denotes the subspace spanned by $S$.

\begin{theorem}[Huang, Tsai \cite{Huang21}]\label{thm: two maps preserving trace}
Let  $\phi: V_1\to W_1$ and $\psi: V_2\to W_2$  be two maps between subsets of matrix spaces over a field $\F$ such that:
\begin{enumerate}
\item $\dim \langle V_1\rangle=\dim \langle V_2\rangle\ge \max\{\dim \langle W_1\rangle, \dim \langle W_2\rangle\}$.
\item %The matrix products
$AB$  are well-defined square matrices for $(A,B)\in (V_1\times V_2)\cup(W_1\times W_2)$.
\item \label{tr invertible}
 If $A\in \langle V_1\rangle $ satisfies that $\tr(AB)=0$ for all $B\in \langle V_2\rangle$, then $A=0$.
\item $\phi$ and $\psi$ satisfy that
\begin{equation}
\tr (\phi(A)\psi(B))=\tr(AB),\quad A\in V_1,\ B\in V_2.
\end{equation}
\end{enumerate}
Then
%\begin{equation}
$\dim \langle V_1\rangle =\dim \langle V_2\rangle =\dim \langle W_1\rangle =\dim \langle W_2\rangle $
%\end{equation}
and $\phi$ and $\psi$ can be extended to bijective linear map $\widetilde\phi : \langle V_1\rangle \to  \langle W_1\rangle$ and $\widetilde\psi: \langle V_2\rangle \to  \langle W_2\rangle$, respectively, such that
\begin{equation}\label{extended map preserve trace}
\tr (\wt{\phi}(A)\wt{\psi}(B))=\tr(AB),\quad   A\in \langle V_1\rangle,\ B\in \langle V_2\rangle.
\end{equation}
%% Moreover, if $A_1,\ldots,A_m$ (resp. $B_1,\ldots, B_m$) is a basis of $\langle V_1\rangle$ (resp. $\langle V_2\rangle$), such that $\tr(A_iB_j)=\delta_{i,j}$, then $\widetilde\phi(A_1),\ldots,\widetilde\phi(A_m)$ (resp. $\widetilde\psi(B_1),\ldots,\widetilde\psi(B_m)$) is a basis of $\langle W_1\rangle$ (resp. $\langle W_2\rangle$) such that $\tr(\widetilde\phi(A_i)\widetilde\psi(B_j))=\delta_{i,j}$.
\end{theorem}

A subset $V$ of $\mathcal{M}_n$   is closed under conjugate transpose if $\{A^*: A\in V\}\subseteq V$.
A real or complex  matrix space $V$ is closed under conjugate transpose if and only if $V$ equals the direct sum of
its subspace of Hermitian  matrices and its subspace of skew-Hermitian  matrices.

\begin{corollary}[Huang, Tsai \cite{Huang21}]\label{thm: trace preserver on special spaces}
Let $V$ be a subset of $\mathcal{M}_n$ closed under conjugate transpose.
Suppose two maps $\phi,\psi:V\to V$   satisfy that
\begin{equation}\label{general two maps preserve trace}
\tr(\phi(A)\psi(B))=\tr(AB),\quad  A, B\in V.
\end{equation}
Then $\phi$ and $\psi$   can be extended to linear bijections on $\langle V\rangle$.
Moreover, when $V$ is a vector space,
 every linear bijection $\phi:V\to V$ corresponds to a unique  linear bijection $\psi:V\to V$ such that \eqref{general two maps preserve trace} holds.
Explicitly, given   an orthonormal basis $\{A_1,\ldots,A_\ell\}$ of $V$ with respect to the inner product $\langle A,B\rangle =\tr(A^*B)$,
$\psi$ is defined by
 $\psi(A_i) =B_i$  in which $\{ B_1,\ldots, B_\ell\}$
 is a basis of $V$ with $\tr(\phi(A_i^*) B_j)=\delta_{i,j}$ for all $i,j\in\{1,\ldots,\ell\}$.
\end{corollary}

%% Furthermore, if $\{A_1,\ldots,A_m\}$ is an orthonormal basis on $V$ with respect to the inner product $\langle A,B\rangle =\tr(A^*B)$, then $\phi$ and $\psi$  satisfy \eqref{general two maps preserve trace} if and only if they can be determined by  a pair of   bases  $\{\phi(A_1^*),\ldots,\phi(A_m^*)\}$ and $\{\psi(A_1),\ldots,\psi(A_m)\}$ on $V$  that satisfy $\tr(\phi(A_i^*)\psi(A_j))=\delta_{i,j}$ for all $i,j\in\{1,\ldots,m\}$.

Corollary \ref{thm: trace preserver on special spaces} shows that when a matrix space $V$ is closed under conjugate transpose,  every  linear bijection $\phi:V\to V$  corresponds to a unique
 linear bijection $\psi:V\to V$ that makes \eqref{general two maps preserve trace} hold.
 The next natural thing is to determine $\phi$ and $\psi$ that satisfy 
$\tr(\phi(A)\psi(B)^k)=\tr(AB^k)$ for  a fixed $k\in\Z\setminus\{0,1\}.$

From now on, we focus on the fields  $\F=\C$ or $\R$. 
%Many results can be easily extended to other subfields of $\C$. 

\section{$k$-power linear preservers and trace of power-product preservers on $\mathcal{M}_n$ and $\mathcal{M}_n(\R)$}

\subsection{$k$-power preservers on $\mathcal{M}_n$ and $\mathcal{M}_n(\R)$}

Chan and Lim described the linear $k$-power preservers on $\mathcal{M}_n$ and $\mathcal{M}_n(\R)$ for $k\ge 2$ in  \cite[Theorem 1]{ChanLim}  as follows.

%Let char $F=0$ or char $F>k$

\begin{theorem} \label{thm: Mn ChanLim} (Chan, Lim  \cite{ChanLim}) Let an integer $k\geq2$. Let $\F$ be a field with $\operatorname{char}(\F)=0$ or
$\operatorname{char}(\F)>k$. Suppose that $\psi: \mathcal{M}_n(\F)\to \mathcal{M}_n(\F)$ is a nonzero linear operator such that $\psi(A^k)=\psi(A)^k$ for all $A\in \mathcal{M}_n(\F)$. Then there exist   $\lambda\in\F$ with $\lambda^{k-1}=1$ and an invertible matrix $P\in \mathcal{M}_n(\F)$ such that
\be{eqnarray}
{ \label{linear bijection preserving power 1}
\psi(A)=\lambda PAP^{-1}, && A\in \mathcal{M}_n(\F),\quad\text{or}
\\ \label{linear bijection preserving power 2}
\psi(A)=\lambda PA^tP^{-1}, &&  A\in \mathcal{M}_n(\F).
}
\end{theorem}

 \eqref{linear bijection preserving power 1} and \eqref{linear bijection preserving power 2} need not hold if $\psi$ is zero or is 
 a map on a subspace of $\mathcal{M}_n(\F)$. The following are two examples. Another example can be found in maps on $\mathcal{D}_n(\F)$ (Theorem \ref{thm: Dn k power}). 

\be{example}
{
The zero map $\psi(A)\equiv 0$    clearly satisfies $\psi(A^k)=\psi(A)^k$ for all $A\in \mathcal{M}_n$ but they are not of the form  \eqref{linear bijection preserving power 1} or \eqref{linear bijection preserving power 2}.
}

\be{example}
{
Let $n=k+m$, $k,m\ge 2$, and consider the operator $\psi$ on the subspace $W=\mathcal{M}_{k}\oplus \mathcal{M}_{m}$ of $\mathcal{M}_n$ defined by
$\psi(A\oplus B)=A\oplus B^t$ for $A\in \mathcal{M}_{k}$ and $B\in\mathcal{M}_{m}.$
Then  $\psi(A^k)=\psi(A)^k$ for all $A\in W$ and $k\in\Z_+$, but $\psi$ is not of the form \eqref{linear bijection preserving power 1} or \eqref{linear bijection preserving power 2}.
}

 We now generalize Theorem \ref{thm: Mn ChanLim} 
to include negative integers $k$ and to assume the   $k$-power preserving condition $\psi(A^k)=\psi(A)^k$ only  on matrices nearby the identity.

\begin{theorem} \label{thm: Mn k-power negative}
Let $\F=\C$ or $\R$. Let an integer $k\in\Z\setminus\{0,1\}$. 
Suppose that $\psi: \mathcal{M}_n(\F)\to \mathcal{M}_n(\F)$ is a nonzero linear map such that $\psi(A^k)=\psi(A)^k$ for all $A$ in an open neighborhood of $I_n$ consisting of invertible matrices. Then there exist   $\lambda\in\F$ with $\lambda^{k-1}=1$ and an invertible matrix $P\in \mathcal{M}_n(\F)$ such that
\be{eqnarray}
{\label{Mn preserving k-power 1}
\psi(A)=\lambda PAP^{-1}, &&   A\in \mathcal{M}_n(\F),\quad\text{or}
\\ \label{Mn preserving k-power 2}
\psi(A)=\lambda PA^tP^{-1}, &&   A\in \mathcal{M}_n(\F).
}
\end{theorem}

\begin{proof}  We prove for the case $\F=\C$. The   case $\F=\R$ can be done similarly.
Obviously, $\psi(I_n)=\psi(I_n^k)=\psi(I_n)^k$.
\begin{enumerate}
\item First suppose $k\ge 2$.
For each $A\in \mathcal{M}_n$, there exists $\epsilon>0$ such that for all $x\in\C$ with $|x|<\epsilon$, the following two power series converge and equal:
\be{eqnarray}
{
(\psi(I_n+xA))^k &=&
 \psi(I_n)+x\left(\sum_{i=0}^{k-1}\psi(I_n)^{i}\psi(A)\psi(I_n)^{k-1-i}\right)   \notag
 \\	&&\quad +x^2\left(\sum_{i=0}^{k-2}\sum_{j=0}^{k-2-i}\psi(I_n)^{i}\psi(A)\psi(I_n)^{j}\psi(A)\psi(I_n)^{k-2-i-j}\right)+\cdots \qquad \label{Mn k-power I+xA 1}
\\
\psi((I_n+xA)^k)
&=&
\psi(I_n)+xk\psi(A)+x^2\frac{k(k-1)}{2} \psi(A^2)+\cdots   \label{Mn k-power I+xA 2}
}
Equating  degree one terms above, we get
\be{equation}{\label{Mn: k psi(A)}
k\psi(A)=\sum_{i=0}^{k-1}\psi(I_n)^{i}\psi(A)\psi(I_n)^{k-1-i}.
}
Applying \eqref{Mn: k psi(A)}, we have
\be{equation}{
k\psi(I_n)\psi(A)-k\psi(A)\psi(I_n)=
\psi(I_n)^{k}\psi(A)-\psi(A)\psi(I_n)^{k}
=\psi(I_n)\psi(A)-\psi(A)\psi(I_n).
}
Hence $\psi(I_n)\psi(A)=\psi(A)\psi(I_n)$ for $A\in \mathcal{M}_n$, that is, $\psi(I_n)$ commutes with the range of $\psi$.

Now equating degree two terms of \eqref{Mn k-power I+xA 1} and \eqref{Mn k-power I+xA 2} and taking into account that $k\not\in\{0,1\}$, we have
\be{equation}{\label{Mn psi^2 equality}
\psi(I_n)^{k-2}\psi(A)^2=\psi(A^2).
}
Define $\psi_1(A)=\psi(I_n)^{k-2}\psi(A)$ for $A\in \mathcal{M}_n$. Then $\psi_1(A^2)=(\psi_1(A))^2$ for all $A\in \mathcal{M}_n$.
\eqref{Mn: k psi(A)} and the assumption that $\psi$ is nonzero imply that $\psi(I_n)\ne 0$. So
$\psi_1(I_n)\psi(I_n)=\psi(I_n)^k=\psi(I_n)\ne 0$. Thus $\psi_1(I_n)\ne 0$ and $\psi_1$ is nonzero.
By  Theorem \ref{thm: Mn ChanLim}, there exists   an invertible $P\in \mathcal{M}_n$ such that
$\psi_1(A)= PAP^{-1}$ for $A\in \mathcal{M}_n$ or $\psi_1(A)=  PA^t P^{-1}$ for $A\in \mathcal{M}_n$.
Moreover, $\psi(I_n)$ commutes with all $\psi_1(A)$, so that $\psi(I_n)=\lambda I_n$ for a $\lambda\in\C$.
By $I_n=\psi_1(I_n)=\psi(I_n)^{k-1}$, we get $\lambda^{k-1}=1$. Therefore, $\psi(A)=\lambda \psi_1(A)$.
We get \eqref{Mn preserving k-power 1} and \eqref{Mn preserving k-power 2}.

\item Next Suppose $k<0$. For every $A\in \mathcal{M}_n$, the power series expansions of
$(\psi(I_n+xA))^{-k}$ and $\psi((I_n+xA)^k)^{-1}$ are equal when $|x|$ is sufficiently small:
\be{eqnarray}
{
(\psi(I_n+xA))^{-k} &=&
 \psi(I_n)^{-1}+x\left(\sum_{i=0}^{-k-1}\psi(I_n)^{i}\psi(A)\psi(I_n)^{-k-1-i}\right)   +\cdots
 \qquad \label{Mn k-power I+xA 3}
\\
\psi((I_n+xA)^k)^{-1}
&=&
\psi(I_n)^{-1}-xk\psi(I_n)^{-1}\psi(A)\psi(I_n)^{-1} +\cdots   \label{Mn k-power I+xA 4}
}
Equating  degree one terms of \eqref{Mn k-power I+xA 3} and \eqref{Mn k-power I+xA 4}, we get
\be{equation}{
-k\psi(I_n)^{-1}\psi(A)\psi(I_n)^{-1}=\sum_{i=0}^{-k-1}\psi(I_n)^{i}\psi(A)\psi(I_n)^{-k-1-i}.
}
Therefore,
\be{eqnarray}{
&&-k(\psi(A)\psi(I_n)^{-1}-\psi(I_n)^{-1}\psi(A))
\notag\\
&=& \psi(I_n)\left(\sum_{i=0}^{-k-1}\psi(I_n)^{i}\psi(A)\psi(I_n)^{-k-1-i}\right)-\left(\sum_{i=0}^{-k-1}\psi(I_n)^{i}\psi(A)\psi(I_n)^{-k-1-i}\right)\psi(I_n)
\notag\\
&=& \psi(I_n)^{-k}\psi(A)-\psi(A)\psi(I_n)^{-k}= \psi(I_n)^{-1}\psi(A)-\psi(A)\psi(I_n)^{-1}.
}
We get $\psi(I_n)^{-1}\psi(A)=\psi(A)\psi(I_n)^{-1}$ for $A\in \mathcal{M}_n$. So $\psi(I_n)^{-1}$ and $\psi(I_n)$ commute with the range of $\psi$.
The following power series are equal for every $A\in \mathcal{M}_n$ when $|x|$ is sufficiently small:
\be{eqnarray}
{
(\psi(I_n+xA))^{k} &=&
 \psi(I_n) +x k\psi(I_n)^{k-1}\psi(A) +x^2\frac{k(k-1)}{2}\psi(I_n)^{k-2}\psi(A)^2   +\cdots
 \quad \label{Mn k-power I+xA 5}
\\
\psi((I_n+xA)^k)
&=&
\psi(I_n) +xk \psi(A)  +x^2\frac{k(k-1)}{2}\psi(A^2)+\cdots   \label{Mn k-power I+xA 6}
}
Equating  degree two terms of \eqref{Mn k-power I+xA 5} and \eqref{Mn k-power I+xA 6}, we get
$\psi(I_n)^{k-2}\psi(A)^2=\psi(A^2).$ Let $\psi_1(A):=\psi(I_n)^{k-2}\psi(A)=\psi(I_n)^{-1}\psi(A)$.
Then $\psi_1(A)^2=\psi_1(A^2)$ and $\psi_1$ is nonzero. Using  Theorem \ref{thm: Mn ChanLim}, we can get \eqref{Mn preserving k-power 1} and \eqref{Mn preserving k-power 2}.
\qedhere
\end{enumerate}
\end{proof}

\subsection{Trace of power-product preserers on $\mathcal{M}_n$ and $\mathcal{M}_n(\R)$}

 Corollary \ref{thm: trace preserver on special spaces}  shows that every linear bijection $\phi:\mathcal{M}_n(\F)\to \mathcal{M}_n(\F)$ corresponds to another linear bijection $\psi:\mathcal{M}_n(\F)\to \mathcal{M}_n(\F)$ such that
$\tr(\phi(A)\psi(B))=\tr(AB)$ for all $A, B\in \mathcal{M}_n(\F)$.
When $m\ge 3$,   maps $\phi_1,\cdots,\phi_m$ on $\mathcal{M}_n(\F)$  that
satisfy $\tr(\phi_1(A_1)\cdots\phi_m(A_m))=\tr(A_1\cdots A_m)$ 
for $A_1,\ldots, A_m\in \mathcal{M}_n(\F)$  are determined in \cite{Huang21}.

If two  maps on $\mathcal{M}_n(\F)$ satisfy  the following  trace condition about $k$-powers, then they have specific forms.

\begin{theorem}\label{Two traces}
Let $\F=\C$ or $\R$.
Let  $k \in\Z\setminus\{0,1\}$. Let $S$ be an open neighborhood of $I_n$ consisting of invertible matrices.
Then two maps $\phi, \psi: \mathcal{M}_n(\F)\to \mathcal{M}_n(\F)$ satisfy that
\be{equation}
{\label{Mn trace two maps power}
 \tr (\phi(A)\psi(B)^k)=\tr (AB^k),%\quad  A \in \mathcal{M}_n(\F),   B\in S,
}
\begin{enumerate}
\item for all $A\in \mathcal{M}_n(\F),\ B\in S$, and $\psi$ is linear, or
\item for all $A, B\in S$ and both $\phi$ and $\psi$ are  linear,
\end{enumerate}
if and only if  $\phi$ and $\psi$ take the following forms:
\begin{enumerate}
\item[(a)] When $k=-1$, there exist invertible matrices $P, Q\in \mathcal{M}_n(\F)$ such that
\be{equation}{\label{Mn two map preserve traces k=-1}
\be{cases}{
\phi(A) =   PAQ\\
\psi(B) =  PBQ
}\quad \text{or}\quad
\be{cases}{
\phi(A) =   PA^t Q\\
\psi(B) =  PB^t Q
}
\quad  A, B\in \mathcal{M}_n(\F).
}

\item[(b)]
When $k\in\Z\setminus\{-1,0,1\}$, there exist $c\in\F\setminus\{0\}$ and an invertible   matrix $P\in \mathcal{M}_n(\F)$ such that
\be{equation}{\label{Mn two map preserve traces}
\be{cases}{
\phi(A) = c^{-k} PAP^{-1}\\
\psi(B) = cPBP^{-1}
}  \text{ or }
\be{cases}{
\phi(A) = c^{-k} PA^t P^{-1}\\
\psi(B) = cPB^t P^{-1}
}
\quad  A, B\in \mathcal{M}_n(\F).
}

\end{enumerate}
\end{theorem}

\begin{proof}
We prove the case $\F=\C$; the case $\F=\R$ can be done similarly.

Suppose assumption (2) holds. Then for every $A\in \mathcal{M}_n(\F)$, there exists $c\in\F\setminus\{0\}$ such that $I_n-cA\in S$, so that for all $B\in S$:
\be{equation}{
\tr(B^k)= \tr((\phi(I_n-cA)+c\phi(A))\psi(B)^k)
= \tr((I_n-cA)B^k)+c\tr(\phi(A)\psi(B)^k).
}
Thus $\tr(\phi(A)\psi(B)^k)=\tr(AB^k)$  for $A \in  \mathcal{M}_n(\F)$ and $B\in S$, which leads to assumption (1). 

Now we prove the theorem under assumption (1), that is, 
\eqref{Mn trace two maps power} holds for  all $A\in \mathcal{M}_n(\F)$ and $B\in S$, and $\psi$ is linear.
Only the necessary part is needed to prove.

%% When $k$ is odd, define $\wt{\psi}:\mathcal{H}_n\to \mathcal{H}_n$ such that $\wt{\psi}(B)=\psi(B^{1/k})^{k}$. When $k$ is even,

Let $S'=\{B\in \ol{\mathcal{P}_n}: B^{1/k}\in S\}$, which is an open neighborhood of $I_n$ in $\ol{\mathcal{P}_n}$.
Define $\wt{\psi}:S'\to \mathcal{M}_n$ such that $\wt{\psi}(B)=\psi(B^{1/k})^{k}$.
Then \eqref{Mn trace two maps power} implies that
\be{equation}{
\tr(\phi(A)\wt{\psi}(B))=\tr(\phi(A)\psi(B^{1/k})^{k})=\tr (AB),\quad  A\in \mathcal{M}_n,\ B\in S'.
}
The complex span of $S'$ is $\mathcal{M}_n$.
By Theorem \ref{thm: two maps preserving trace},  $\phi$ is bijective linear,
and $\wt{\psi}$  can be extended to a  linear bijection on $\mathcal{M}_n$.  %%In particular, $\wt{\psi}$ maps every basis of $\mathcal{M}_n$ in $\ol{\mathcal{P}_n}$ to a basis of $\mathcal{M}_n$.

%%For convenience, denote $A'=\phi(A)$ and $A''=\psi(A)$ for $A\in \mathcal{H}_n$ in this proof.
The linearity of $\psi$ and \eqref{Mn trace two maps power} imply that for every $B\in \mathcal{M}_n$, there exists $\epsilon>0$ such that
$I_n+xB\in S$ and the power series of $(I_n+xB)^k$ converges whenever $|x|<\epsilon$. Then
\be{equation}{\label{Mn-I+xB}
\tr \left(\phi(A)(\psi(I_n)+x\psi(B))^k\right)=\tr (A(I_n+xB)^k),\quad  A\in \mathcal{M}_n,\  |x|<\epsilon.
}

\begin{enumerate}
\item First suppose $k\ge 2$.
Equating degree one terms and degree $(k-1)$ terms  on both sides of \eqref{Mn-I+xB} respectively, we get the following identities for $A, B\in \mathcal{M}_n$:
\be{eqnarray}{\label{Mn-(k,1)}
  \tr \left(\phi(A)\left (\sum_{i=0}^{k-1}\psi(I_n)^{k-1-i}\psi(B)\psi(I_n)^{i}\right )\right) &=& \tr (kAB),
\\ \label{Mn-(k,k-1)}
 \tr \left (\phi(A)\left (\sum_{i=0}^{k-1} \psi(B)^{i}\psi(I_n)\psi(B)^{k-1-i}\right )\right) &=& \tr (kAB^{k-1}).
}

Let $\{C_i:i=1,\ldots,n^2\}$ be a basis of projection matrices (i.e. $C_i^2=C_i$) in $\mathcal{M}_n$. For example, we may choose the following basis of rank 1 projections:
\be{equation}{\label{Mn projection basis}
\{E_{ii}: 1\le i\le n\}\cup \left\{\frac{1}{\sqrt{2}}(E_{ii}+E_{jj}+\delta E_{ij}+\bar\delta E_{ji}):1\le i<j\le n, \delta\in\{1,\i\}\right\}.
}
By \eqref{Mn trace two maps power} and \eqref{Mn-(k,k-1)}, for $A\in \mathcal{M}_n$ and $i=1,\ldots,n^2$,
\be{equation}{
\tr(k\phi(A)\psi(C_i)^{k})= \tr (kAC_i) =  \tr \left (\phi(A)\left (\sum_{j=0}^{k-1} \psi(C_i)^{j}\psi(I_n)\psi(C_i)^{k-1-j}\right )\right).
}
By the bijectivity of $\phi$,
\be{equation}{
k\psi(C_i)^{k}=\sum_{j=0}^{k-1} \psi(C_i)^{j}\psi(I)\psi(C_i)^{k-1-j}.
}
Therefore, for $i=1,\ldots,n^2$,
\be{eqnarray}{\notag
0&=& \psi(C_i)\left (\sum_{j=0}^{k-1} \psi(C_i)^{j}\psi(I_n)\psi(C_i)^{k-1-j}\right )-\left (\sum_{j=0}^{k-1} \psi(C_i)^{j}\psi(I_n)\psi(C_i)^{k-1-j}\right )\psi(C_i)
\\\label{Mn-psi(I)psi(C_i)^k}
&=&
\psi(C_i)^{k}\psi(I_n)-\psi(I_n)\psi(C_i)^{k}.
}
Since
\be{equation}{
\tr(A\psi(C_i)^k)=\tr(\phi^{-1}(A)C_i^k)=\tr(\phi^{-1}(A)C_i),\quad  A\in \mathcal{M}_n,\ i=1,\ldots,n^2,
}
the only matrix $A\in \mathcal{M}_n$ such that  $\tr(A\psi(C_i)^k)=0$ for all $i\in\{1,\ldots,n^2\}$ is the zero matrix.
So $\{\psi(C_i)^{k}:i=1,\ldots,n^2\}$ is a basis of $\mathcal{M}_n$.  \eqref{Mn-psi(I)psi(C_i)^k} implies that  $\psi(I_n)=cI_n$ for certain $c\in\C\setminus\{0\}$.

 \eqref{Mn-(k,1)}  shows that
\be{equation}{
c^{k-1}\tr(\phi(A)\psi(B))=\tr(AB),\quad A, B\in \mathcal{M}_n.
}
Therefore,
\be{equation}{
c^{k-1}\tr(\phi(A)\psi(B^k))=\tr(AB^k)=\tr(\phi(A)\psi(B)^k),\quad A\in \mathcal{M}_n,\ B\in S.
}
The bijectivity of $\phi$ shows that $c^{k-1}\psi(B^k)=\psi(B)^k$ for $B\in S$, that is,
\be{equation}{\label{Mn c^-1 psi}
c^{-1}\psi(B^k)=[c^{-1}\psi(B)]^k,\quad B\in S.
}
Notice that $c^{-1}\psi(I_n)=I_n$.  By  Theorem \ref{thm: Mn k-power negative},
 there is an invertible $P\in \mathcal{M}_n$ such that $\psi$ is of the form $\psi(B) = cPBP^{-1}$ or $\psi(B) = cPB^t P^{-1}$ for $B\in \mathcal{M}_n$.
Consequently, we get \eqref{Mn two map preserve traces}.

\item Now suppose $k<0$. Then $\psi(I_n)$ is invertible. For every $B\in \mathcal{M}_n$ and sufficiently small $x$, we have the power series expansion:
\be{eqnarray}{\notag
&& (\psi(I_n)+x\psi(B))^k
\\ \notag
&=& \left[(I_n+x\psi(I_n)^{-1}\psi(B))^{-1}\psi(I_n)^{-1}\right]^{|k|}
\\ \notag
&=& \left[(I_n-x\psi(I_n)^{-1}\psi(B)+x^2\psi(I_n)^{-1}\psi(B)\psi(I_n)^{-1}\psi(B)+\cdots) \psi(I_n)^{-1}\right]^{|k|}
\\ \notag
&=& \psi(I_n)^{k}-x\left(\sum_{i=1}^{|k|} \psi(I_n)^{-i}\psi(B)\psi(I_n)^{k-1+i}\right)
\\ \label{Mn psi(I+xB)}
&&\quad
+x^2\left(\sum_{i=1}^{|k|}\sum_{j=1}^{|k|+1-i}\psi(I_n)^{-i}\psi(B)\psi(I_n)^{-j}\psi(B)\psi(I_n)^{k-2+i+j}  \right)+\cdots
}
Equating degree one terms and degree two terms of \eqref{Mn-I+xB} respectively and using \eqref{Mn psi(I+xB)}, we get the following identities for $A, B\in \mathcal{M}_n$:
\be{eqnarray}{\label{Mn negative power deg 1}
 \tr\left (\phi(A)\left(\sum_{i=1}^{|k|} \psi(I_n)^{-i}\psi(B)\psi(I_n)^{k-1+i}\right)\right) &=& \tr (|k|AB),
\\ \label{Mn negative power deg 2}
  \tr\left (\phi(A)\left(\sum_{i=1}^{|k|}\sum_{j=1}^{|k|+1-i}\psi(I_n)^{-i}\psi(B)\psi(I_n)^{-j}\psi(B)\psi(I_n)^{k-2+i+j}  \right) \right) &=& \tr \left(\frac{k(k-1)}{2} AB^2 \right).\qquad
}
\eqref{Mn negative power deg 1} and \eqref{Mn trace two maps power} imply that
\be{equation}{\label{Mn power identity}
 \sum_{i=1}^{|k|} \psi(I_n)^{-i}\psi(B^k)\psi(I_n)^{k-1+i} =|k|\psi(B)^k,\quad  B\in S.
}

Let $F_r(B)$ denote the degree $r$ coefficient in the power series of $(\psi(I_n)+x\psi(B))^k$. Then \eqref{Mn negative power deg 1} and \eqref{Mn negative power deg 2}
show that:
\be{equation}{\label{Mn T2 T1}
\frac{k-1}{2} F_1(B^2)=F_2(B),\quad  B\in \mathcal{M}_n.
}

Denote $\psi_1(B):=\psi(B)\psi(I_n)^{-1}$. We discuss the cases $k=-1$ and $k\ne -1$.

\begin{enumerate}
\item
When $k=-1$, \eqref{Mn T2 T1} leads to
\be{equation}{\label{Mn psi B^2 identity}
\psi(B^2)=\psi(B)\psi(I_n)^{-1}\psi(B),\quad  B\in \mathcal{M}_n.
 }
So  $\psi_1(B^2)= \psi_1(B)^2$ for $B\in \mathcal{M}_n$. Note that $\psi_1(I_n)=I_n$.  By Theorem \ref{thm: Mn k-power negative}, there exists an invertible $P\in \mathcal{M}_n$ such that
$\psi_1(B)=PBP^{-1}$ or $\psi_1(B)= PB^t P^{-1}$ for $B\in \mathcal{M}_n$.
Let $Q:=P^{-1}\psi(I_n)$. Then $Q$ is invertible, and $\psi(B)=PBQ$ or $\psi(B)=PB^t Q$ for $B\in \mathcal{M}_n$.
Using \eqref{Mn trace two maps power}, we get \eqref{Mn two map preserve traces k=-1}.

\item Suppose the integer $k<-1$.
Then \eqref{Mn T2 T1} implies that
\be{equation}{
\frac{k-1}{2} (\psi(I_n)^{-1}F_1(B^2)-F_2(B^2)\psi(I_n)^{-1})=\psi(I_n)^{-1}F_2(B)-F_2(B)\psi(I_n)^{-1},
}
which gives
\be{eqnarray}{
&& \frac{1-k}{2}\left(\psi(I_n)^{k}\psi(B^2)-\psi(B^2)\psi(I_n)^{k}\right)
\notag \\
&=& \psi(I_n)^{k}\psi(B)\psi(I_n)^{-1}\psi(B)-\psi(B)\psi(I_n)^{-1}\psi(B)\psi(I_n)^{k}.
\label{Mn psi(I)B^2}
}
In other words, for $B\in \mathcal{M}_n$:
\be{equation}{\label{Mn B^2 equality}
\psi(I_n)^{k}\left(\frac{1-k}{2}\psi_1(B^2)-\psi_1(B)^2\right)=\left(\frac{1-k}{2}\psi_1(B^2)-\psi_1(B)^2\right)\psi(I_n)^{k}.
}
Let $B=I_n+xE$ for an arbitrary matrix $E\in \mathcal{M}_n$. Then \eqref{Mn B^2 equality} becomes
\be{eqnarray}{\label{Mn B^2 expansions}
&&\psi(I_n)^{k}\left[x(-1-k)\psi_1(E)+x^2\left(\frac{1-k}{2}\psi_1(E^2)-\psi_1(E)^2\right)\right]
\notag \\
&=&
\left[x(-1-k)\psi_1(E)+x^2\left(\frac{1-k}{2}\psi_1(E^2)-\psi_1(E)^2\right)\right]\psi(I_n)^{k}.
}
The equality on degree one terms shows that $\psi(I_n)^{k}$ commutes with all $\psi_1(E)$. Hence $\psi(I_n)^k$ commutes with the range of $\psi$.
  \eqref{Mn negative power deg 1} can be rewritten as
\be{equation}{
\tr\left (\left(\sum_{i=1}^{|k|}  \psi(I_n)^{k-1+i}\phi(A)\psi(I_n)^{-i} \right)  \psi(B)\right) = \tr (|k|AB),\quad  A, B\in \mathcal{M}_n.
}
By Theorem \ref{thm: two maps preserving trace}, $\psi$ is a linear bijection  and its range is $\mathcal{M}_n$. So
$\psi(I_n)^k=\mu I_n$ for certain $\mu\in\C$.

Now by \eqref{Mn power identity}, for $B\in S$:
\be{eqnarray}{
&&|k|\psi(I_n)\psi(B)^k - |k|\psi(B)^k\psi(I_n)
\notag\\
&=& \psi(I_n)\left( \sum_{i=1}^{|k|} \psi(I_n)^{-i}\psi(B^k)\psi(I_n)^{k-1+i}\right)-\left( \sum_{i=1}^{|k|} \psi(I_n)^{-i}\psi(B^k)\psi(I_n)^{k-1+i}\right)\psi(I_n)
\notag\\
&=&
\psi(B^k)\psi(I_n)^k-\psi(I_n)^k\psi(B^k)=0
}
So $\psi(I_n)$ commutes with $\psi(B)^k$ for $B\in S$.
In particular, $\psi(I_n)$ commutes with  $\wt{\psi}(B)=\psi(B^{1/k})^k$ for $B\in S'$.  The complex span of $S'$ is $\mathcal{M}_n$, and $\wt{\psi}$ can be extended to a linear bijection on $\mathcal{M}_n$. Hence $\psi(I_n)=c I_n$ for certain $c \in \C\setminus\{0\}$.
By \eqref{Mn power identity}, we get $\psi_1(B^k)=\psi_1(B)^k$ for $B\in S$. Note that $\psi_1(I_n)=I_n$.
By Theorem \ref{thm: Mn k-power negative}, there is an invertible $P\in \mathcal{M}_n$ such that
$\psi_1(B)= PBP^{-1}$ or $\psi_1(B)= PB^t P^{-1}$.
Then $\psi (B)=c PBP^{-1}$ or $\psi (B)=c PB^t P^{-1}$.
Using \eqref{Mn trace two maps power}, we get \eqref{Mn two map preserve traces}.
\qedhere

\end{enumerate}
\end{enumerate}
\end{proof}

\begin{remark}
The following modifications could be applied to the proof of  Theorem \ref{Two traces} for $\F=\R$:
\begin{enumerate}
\item
Let $S'$ be the collection of matrices $A=QDQ^{-1}$, in which $D$ is nonnegative diagonal and $Q\in \mathcal{M}_n(\R)$ is invertible, such that $A^{1/k}=QD^{1/k}Q^{-1}\in S$.

\item
We may choose the following basis of rank 1 projections of $\mathcal{M}_n(\R)$ to substitute
\eqref{Mn projection basis}:
\be{multline}{
\{E_{ii}: 1\le i\le n\}\cup \left\{\frac{1}{\sqrt{2}}(E_{ii}+E_{jj}+E_{ij}+E_{ji}):1\le i<j\le n\right\}
\\
\cup  \left\{\omega_1 E_{ii}+\omega_2 E_{jj}+E_{ij}-E_{ji}:1\le i<j\le n\right\},
}
in which $\omega_1, \omega_2$ are distinct roots of $x^2-x-1=0.$
\end{enumerate}
\end{remark}

The  arguments in the above proof will be applied analogously to maps on the other sets discussed in this paper.

%% It is unclear whether or not  \eqref{Mn two map preserve traces} still holds when the assumption ``$\psi: \mathcal{M}_n\to \mathcal{M}_n$ is linear''  in Theorem \ref{Two traces} is removed.  However, the linearity of $\psi$ must be assumed for a similar result on $\mathcal{H}_n$. See Remark \ref{rem: Hn power linearity}.

%%\section{Maps $\mathcal{H}_n\to \mathcal{H}_n$ and $\mathcal{P}_n\to \ol{\mathcal{P}_n}$ preserving trace of products}

\section{$k$-power linear preservers and trace of power-product preservers on $\mathcal{H}_n$}

%% The following result for maps $\mathcal{H}_n\to \mathcal{H}_n$ is similar to that of $\mathcal{P}_n\to \mathcal{P}_n$ for $m\ge 3$.

\subsection{$k$-power linear preservers on $\mathcal{H}_n$}

We give a result that determine linear operators on $\mathcal{H}_n$ that satisfy $\psi(A^k)=\psi(A)^k$ on a neighborhood of $I_n$ in  $\mathcal{H}_n$  for certain $k\in\Z\setminus\{0,1\}$.

\be{theorem}
{\label{thm: Hn linear map preserving k power}
Fix $k\in\Z\setminus\{0,1\}$.
A nonzero linear map $\psi: \mathcal{H}_n\to \mathcal{H}_n$ satisfies that
\be{equation}{\label{Hn preserve k power}
\psi(A^k)=\psi(A)^k
}
on an open neighborhood of $I_n$ consisting of invertible matrices 
if and only if $\psi$ is of the following forms for certain unitary matrix $U\in \mathcal{M}_n$:
\begin{enumerate}
\item When $k$ is even,
\be{equation}
{\label{Hn: trace power even}
\psi(A)=U^*AU,\quad  A\in \mathcal{H}_n;\quad\text{or}\quad \psi(A)=U^*A^t U,\quad  A\in \mathcal{H}_n.
}
\item When $k$ is odd,
\be{equation}
{\label{Hn: trace power odd}
\psi(A)=\pm U^*AU,\quad  A\in \mathcal{H}_n;\quad\text{or}\quad \psi(A)=\pm U^*A^t U,\quad  A\in \mathcal{H}_n.
}
\end{enumerate}
}

\begin{proof}
It suffices to prove the necessary part. Suppose \eqref{Hn preserve k power} holds on an open neighborhood $S$ of $I_n$ in $\mathcal{H}_n$.

\begin{enumerate}
\item First assume $k\ge 2$.
Replacing $\mathcal{M}_n$ by $\mathcal{H}_n$ in part (1) of the proof of  Theorem \ref{thm: Mn k-power negative} up to \eqref{Mn psi^2 equality},
we can prove that $\psi(I_n)$ commutes with the range of $\psi$, and $\psi_1(A):=\psi(I_n)^{k-2}\psi(A)$ is a nonzero linear map that satisfies
$\psi_1(A^2)=\psi_1(A)^2$ for $A\in \mathcal{H}_n$.

Every matrix in $\mathcal{M}_n$ can be uniquely expressed as $A+\i B$ for $A,B\in \mathcal{H}_n$.
Extend $\psi_1$ to a map $\widetilde\psi:\mathcal{M}_n\to \mathcal{M}_n$ such that
\be{equation}{
\wt{\psi}(A+\i B)=\psi_1(A)+\i\psi_1(B),\quad A,B\in \mathcal{H}_n.
}
It is straightforward to check that $\wt{\psi}$ is a complex linear bijection. Moreover, for  $A, B\in \mathcal{H}_n$,
\be{eqnarray*}{
\psi_1(AB+BA)
&=& \psi_1((A+B)^2)-\psi_1(A^2)-\psi_1(B^2)
\\
&=&  \psi_1(A+B)^2-\psi_1(A)^2-\psi_1(B)^2
\\
&=& \psi_1(A)\psi_1(B)+\psi_1(B)\psi_1(A).
}
It implies that
$$
\wt{\psi}((A+\i B)^2)=\wt{\psi}(A+\i B)^2,\quad   A,B\in \mathcal{H}_n.
$$
By Theorem \ref{thm: Mn ChanLim}, there is an invertible matrix $U\in \mathcal{M}_n$ such that
 \begin{enumerate}
\item
 $\wt{\psi}(A)=UAU^{-1}$  for all $A\in \mathcal{M}_n$, or
\item
 $\wt{\psi}(A)=UA^t U^{-1}$  for all $A\in \mathcal{M}_n$.
\end{enumerate}

First suppose  $\wt{\psi}(A)=UAU^{-1}$.  The restriction of $\wt{\psi}$ on $\mathcal{H}_n$ is $\psi_1:\mathcal{H}_n\to \mathcal{H}_n$. Hence
for $A\in \mathcal{H}_n$, we have $UAU^{-1}=(UAU^{-1})^*=U^{-*}AU^*$; then $U^*UA=AU^*U$ for all $A\in \mathcal{H}_n$, which shows that $U^*U=cI_n$ for certain $c\in\R^+.$
By adjusting a scalar if necessary, we may assume that $U$ is unitary. So $\psi(I_n)^{k-2}\psi(A)=UAU^*$. Then $\psi(I_n)^{k-1}=I_n$, so that
$\psi(I_n)=I_n$ when $k$ is even and $\psi(I_n)\in\{I_n,-I_n\}$ when $k$ is odd. Thus $\psi(A)=UAU^*$  when $k$ is even and
 $\psi(A)=\pm UAU^*$  when $k$ is odd.
Similarly for the case $\wt{\psi}(A)=UA^t U^{-1}$.
Therefore, \eqref{Hn: trace power even} or \eqref{Hn: trace power odd} holds.

\item Now assume that $k<0$. Replacing $\mathcal{M}_n$ by $\mathcal{H}_n$ in part (2)   of the proof of Theorem \ref{thm: Mn k-power negative},
we can show that $\psi(I_n)$ commutes with the range of $\psi$, and furthermore  the nonzero linear map $\psi_1(A):=\psi(I_n)^{-1}\psi(A)$
satisfies that $\psi_1(A^2)=\psi_1(A)^2$.  By arguments in the preceding paragraphs, we get \eqref{Hn: trace power even} or \eqref{Hn: trace power odd}.
\qedhere
\end{enumerate}
\end{proof}

\subsection{Trace of power-product preservers on $\mathcal{H}_n$}

By Corollary \ref{thm: trace preserver on special spaces}, every linear bijection $\phi: \mathcal{H}_n\to \mathcal{H}_n$
corresponds to another linear bijection $\psi: \mathcal{H}_n\to \mathcal{H}_n$ such that
 $\tr(\phi(A)\psi(B))=\tr(AB)$ for all $A, B\in \mathcal{H}_n$.
When $m\ge 3$, linear maps $\phi_1,\cdots,\phi_m: \mathcal{H}_n\to \mathcal{H}_n$ that
satisfy $\tr(\phi_1(A_1)\cdots\phi_m(A_m))=\tr(A_1\cdots A_m)$ are characterized in \cite{Huang21}.

\begin{theorem}\label{thm: Hn Two traces}
Let $k\in\Z\setminus\{0,1\}$. Let $S$ be an open neighborhood of $I_n$ in $\mathcal{H}_n$ consisting of invertible Hermitian matrices.
Then two maps $\phi,\psi: \mathcal{H}_n\to \mathcal{H}_n$   satisfy that
\be{equation}
{\label{Hn trace two maps power}
 \tr (\phi(A)\psi(B)^k)=\tr (AB^k), %\quad  A\in \mathcal{H}_n,\   B\in S,
}
\begin{enumerate}
\item for all $A\in \mathcal{H}_n,\ B\in S$, and $\psi$ is linear, or
\item for all $A, B\in S$ and both $\phi$ and $\psi$ are  linear,
\end{enumerate}
if and only if $\phi$ and $\psi$ take the following forms:
\begin{enumerate}
\item[(a)] When $k=-1$, there exist an invertible matrix $P\in \mathcal{M}_n$ and $c\in\{1,-1\}$ such that
\be{equation}{\label{Hn two map preserve traces k is -1}
\be{cases}{
\phi(A) = cP^*AP\\
\psi(B) = cP^*BP
}   A, B\in \mathcal{H}_n;
 \text{ or }
\be{cases}{
\phi(A) = cP^*A^t P\\
\psi(B) = cP^*B^t P
}
  A, B\in \mathcal{H}_n.
}

\item[(b)] When $k\in\Z\setminus\{-1,0,1\}$, there exist a unitary matrix $U\in \mathcal{M}_n$ and $c\in\R\setminus\{0\}$ such that
\be{equation}{\label{Hn two map preserve traces k is not -1}
\be{cases}{
\phi(A) = c^{-k} U^*AU\\
\psi(B) = cU^*BU
}  A, B\in \mathcal{H}_n; \text{ or }
\be{cases}{
\phi(A) = c^{-k} U^*A^t U\\
\psi(B) = cU^*B^t U
}
  A, B\in \mathcal{H}_n.
}
\end{enumerate}
\end{theorem}

\begin{proof}  
Assumption (2) leads to assumption (1) (cf. the proof of Theorem \ref{Two traces}). 
We prove the theorem under assumption (1). 
It suffices to prove the necessary part.
\begin{enumerate}
\item When $k\ge 2$,  in the part (1) of proof of Theorem \ref{Two traces}, through replacing $\mathcal{M}_n$ by $\mathcal{H}_n$, complex numbers by real numbers, and Theorem \ref{thm: Mn ChanLim} or Theorem \ref{thm: Mn k-power negative}  by Theorem \ref{thm: Hn linear map preserving k power}, we can prove that $(\phi,\psi)$ has the forms in \eqref{Hn two map preserve traces k is not -1}.

\item When $k<0$, in the part (2) of proof of Theorem \ref{Two traces}, through replacing $\mathcal{M}_n$ by $\mathcal{H}_n$ and complex numbers by real numbers,
we can get the corresponding equalities of \eqref{Mn psi(I+xB)} $\sim$ \eqref{Mn T2 T1} on $\mathcal{H}_n$.  The case $k<-1$ can be proved completely analogously with the help of Theorem \ref{thm: Hn linear map preserving k power}.

For the case $k=-1$, the equality corresponding to \eqref{Mn T2 T1} can be simplified as
\be{equation}{\label{Hn T2 T1}
\psi(B^2)=\psi(B)\psi(I_n)^{-1}\psi(B),\quad  B\in \mathcal{H}_n.
}
Let $\psi_1(B):=\psi(I_n)^{-1}\psi(B)$. Then $\psi_1:\mathcal{H}_n\to \mathcal{M}_n$ is a nonzero real linear map
that satisfies $\psi_1(B^2)=\psi_1(B)^2$ for $B\in \mathcal{H}_n$.
 Extend $\psi_1$ to a complex linear map $\wt{\psi}:\mathcal{M}_n\to \mathcal{M}_n$ such that
\be{equation}{
\wt{\psi}(A+\i B):=\psi_1(A)+\i \psi_1(B),\quad  A, B\in \mathcal{H}_n.
}
Similarly to the arguments in part (1) of the proof of Theorem \ref{thm: Hn linear map preserving k power},
we have $\wt{\psi}((A+\i B)^2)=(\wt{\psi}(A+\i B))^2$ for all $A, B\in \mathcal{H}_n$. Using Theorem \ref{thm: Mn k-power negative}
and the fact that $\wt{\psi}(I_n)=\psi_1(I_n)=I_n$, we can prove that there is an invertible $P\in \mathcal{M}_n$ such that  for all $B\in \mathcal{H}_n$, either
$\psi_1(B)=P^{-1}BP$ or $\psi_1(B)=P^{-1} B^t P$. So
\be{eqnarray}{
\psi(B)=\psi(I_n)P^{-1}BP,&& B\in \mathcal{H}_n;\quad \text{or}
\\
\psi(B)=\psi(I_n)P^{-1}B^t P,&& B\in \mathcal{H}_n.
}
If $\psi(B)=\psi(I_n)P^{-1}BP$ for $B\in \mathcal{H}_n$, then
$\psi(I_n)P^{-1}BP = (\psi(I_n)P^{-1}BP)^*=P^*BP^{-*}\psi(I_n),$
which gives
\be{equation}{
(P^{-*}\psi(I_n)P^{-1})B=B(P^{-*}\psi(I_n)P^{-1}),\quad  B\in \mathcal{H}_n.
}
Hence $P^{-*}\psi(I_n)P^{-1}=cI_n$ for certain $c\in\R\setminus\{0\}.$
We have $\psi(I_n)=cP^*P$ so that $\psi(B)=cP^*BP$ for $B\in \mathcal{H}_n$. Similarly for the case $\psi(B)=\psi(I_n)P^{-1}B^t P$.
Adjusting $c$ and $P$ by scalar factors simultaneously, we may assume that $c\in\{1,-1\}$. It implies \eqref{Hn two map preserve traces k is -1}.
\qedhere

%Let $B=I_n+xC$ for $C\in \mathcal{H}_n$ and $x\in\R$ in \eqref{Hn T2 T1} and compare  degree one terms on both sides:
%\be{equation}{2\psi(C)=}

\end{enumerate}

\end{proof}

\begin{remark}\label{rem: Hn power linearity}
Theorem \ref{thm: Hn Two traces} does not hold if  $\psi$ is not assumed to be linear. Let $k$ be a positive even integer.
Let $\wt{\psi}:\mathcal{H}_n\to \mathcal{H}_n$ be any bijective linear map such that $\wt{\psi}(\ol{\mathcal{P}_n})\subseteq \ol{\mathcal{P}_n}$. For example,
$\wt{\psi}$ may be a completely positive map of the form  $\wt{\psi}(B)=\sum_{i=1}^{r} N_i^* B N_i$ for $r\ge 2$, $N_1,\ldots, N_r\in \mathcal{M}_n$ linearly independent, and at least one of $N_1,\ldots, N_r$ is invertible.
By Corollary \ref{thm: trace preserver on special spaces}, there is a linear bijection $\phi:\mathcal{H}_n\to \mathcal{H}_n$ such that
$\tr(\phi(A)\wt{\psi}(B))=\tr(AB)$ for all $A, B\in \mathcal{H}_n$.
Let $\psi:\mathcal{H}_n\to \mathcal{H}_n$ be defined by $\psi(B)=\wt{\psi}(B^k)^{1/k}.$  Then
$$
\tr(\phi(A)\psi(B)^k)=\tr(\phi(A)\wt{\psi}(B^k))=\tr(AB^k),\quad   A, B\in \mathcal{H}_n.
$$
Obviously, $\psi$ may  be non-linear, and the choices of pairs $(\phi,\psi)$ are much more than
those in \eqref{Hn two map preserve traces k is -1} and \eqref{Hn two map preserve traces k is not -1}.
\end{remark}

%
% \textcolor{red}{\bf Note:} We should develop similar results for $\mathcal{M}_n$, $\mathcal{P}_n$, and $\mathcal{S}_n$, etc., and extend to the multiple map cases.

%%%%%%%%%%%%%%%%%%%%%%%%%%%%%%%%%%%%%%

\section{$k$-power linear preservers and trace of power-product preservers on $\mathcal{S}_n$ and $\mathcal{S}_n(\R)$}

\subsection{$k$-power linear preservers on $\mathcal{S}_n$ and $\mathcal{S}_n(\R)$}

Chan and Lim described the linear $k$-power preservers on $\mathcal{S}_n(\F)$ for $k\ge 2$ in  \cite[Theorem 2]{ChanLim}  as follows.

\begin{theorem} \label{thm: Sn ChanLim} (Chan, Lim \cite{ChanLim}) Let an integer $k\geq2$. Let $\F$ be an algebraic closed field with $\operatorname{char}(\F)=0$ or
$\operatorname{char}(\F)>k$.  Suppose that $\psi: \mathcal{S}_n(\F)\to \mathcal{S}_n(\F)$ is a nonzero linear operator such that $\psi(A^k)=\psi(A)^k$ for all $A\in \mathcal{S}_n(\F)$. Then there exist    $\lambda\in\F$ with $\lambda^{k-1}=1$ and an orthogonal matrix $O\in \mathcal{M}_n(\F)$ such that
\be{equation}
{\label{linear bijection preserving power 3}
\psi(A)=\lambda OAO^t, \quad A\in \mathcal{S}_n.
}
\end{theorem}

 We   generalize Theorem \ref{thm: Sn ChanLim}  to include the case $\mathcal{S}_n(\R)$, to include negative integers $k$, and to assume the   $k$-power preserving condition  only  on matrices nearby the identity.

\begin{theorem} \label{thm: Sn k-power negative}
Let   $k\in\Z\setminus\{0,1\}$. Let $\F=\C$ or $\R$.
Suppose that $\psi: \mathcal{S}_n(\F)\to \mathcal{S}_n(\F)$ is a nonzero linear map such that $\psi(A^k)=\psi(A)^k$ for all $A$ in an open neighborhood of $I_n$ in $\mathcal{S}_n(\F)$ consisting of invertible matrices. Then there exist   $\lambda\in\F$ with $\lambda^{k-1}=1$ and an orthogonal matrix $O\in \mathcal{M}_n(\F)$ such that
\be{equation}
{\label{Sn preserving k-power}
\psi(A)=\lambda OAO^{t}, \quad   A\in \mathcal{S}_n(\F).
}
\end{theorem}

\begin{proof}
It suffices to prove the necessary part.
In both $k\ge 2$ and $k<0$ cases, using  analogous arguments as parts (1) and (2) of the proof of Theorem \ref{thm: Mn k-power negative},
 we get that $\psi(I_n)$ commutes with the range of $\psi$,  and the nonzero map
$\psi_1(A):=\psi(I_n)^{k-2}\psi(A)$ satisfies that $\psi_1(A^2)=\psi_1(A)^2$ for $A\in \mathcal{S}_n(\F)$. Then
\be{eqnarray}{
\psi_1(A)\psi_1(B)+\psi_1(B)\psi_1(A)&=& \psi_1(A+B)^2-\psi_1(A)^2-\psi_1(B)^2
\notag \\
&=& \psi_1((A+B)^2)-\psi_1(A^2)-\psi_1(B^2)
\notag \\
&=& \psi_1(AB + BA).
}
In particular, $\psi_1(A)\psi_1(A^r)+\psi_1(A^r)\psi_1(A)=2\psi_1(A^{r+1})$ for $r\in\Z_+$. Using induction, we get
$\psi_1(A^{\ell})=\psi_1(A)^{\ell}$ for all $A\in \mathcal{S}_n(\F)$ and $\ell\in\Z_+$.
By \cite[Corollary 6.5.4]{ZTC}, there is an orthogonal matrix $O\in \mathcal{M}_n(\F)$ such that
$\psi_1(A)=OAO^t$. Since $\psi(I_n)$ commutes with the range of  $\psi_1$, we have $\psi(I_n)=\lambda I_n$ for certain $\lambda\in\F$ in which
$\lambda^{k-1}=1$. So $\psi(A)=\lambda OAO^t$ as in  \eqref{Sn preserving k-power}.
\end{proof}

Obviously, in $\F=\R$ case, \eqref{Sn preserving k-power} has $\lambda=1$ when $k$ is even and $\lambda\in\{1,-1\}$ when $k$ is odd.

\subsection{Trace of power-product preservers on $\mathcal{S}_n$ and $\mathcal{S}_n(\R)$}

 Corollary \ref{thm: trace preserver on special spaces}  shows that every linear bijection $\phi:\mathcal{S}_n(\F)\to \mathcal{S}_n(\F)$ corresponds to another linear bijection $\psi:\mathcal{S}_n(\F)\to \mathcal{S}_n(\F)$ such that
$\tr(\phi(A)\psi(B))=\tr(AB)$ for all $A, B\in \mathcal{S}_n(\F)$.
When $m\ge 3$,   maps $\phi_1,\cdots,\phi_m: \mathcal{S}_n(\F)\to \mathcal{S}_n(\F)$ that
satisfy $\tr(\phi_1(A_1)\cdots\phi_m(A_m))=\tr(A_1\cdots A_m)$ are determined in \cite{Huang21}.

We characterize the trace of power-product preserver for $\mathcal{S}_n(\F)$ here.

\begin{theorem}\label{thm: Sn Two traces}
Let $\F=\C$ or $\R$. Let  $k \in\Z\setminus\{0,1\}$. Let $S$ be an open neighborhood of $I_n$ in $\mathcal{S}_n(\F)$ consisting of invertible matrices.
Then two maps $\phi,\psi: \mathcal{S}_n(\F)\to \mathcal{S}_n(\F)$  satisfy that
\be{equation}
{\label{Sn trace two maps power}
 \tr (\phi(A)\psi(B)^k)=\tr (AB^k), %\quad  A\in \mathcal{S}_n(\F),\   B\in S,
}
\begin{enumerate}
\item for all $A\in \mathcal{S}_n(\F),\ B\in S$, and $\psi$ is linear, or
\item for all $A, B\in S$ and both $\phi$ and $\psi$ are  linear,
\end{enumerate}
if and only if   $\phi$ and $\psi$ take the following forms:
\begin{enumerate}
\item[(a)] When $k=-1$, there exist an invertible matrix $P\in \mathcal{M}_n(\F)$ and $c\in\F\setminus\{0\}$ such that
\be{equation}{\label{Sn two map preserve traces power -1}
\phi(A)=cP A P^t,\quad\psi(B)= cP B P^t,\quad  A, B\in \mathcal{S}_n(\F).
}
We may choose $c=1$ for $\F=\C$ and $c\in\{1,-1\}$ for $\F=\R$.

\item[(b)] When $k\in\Z\setminus\{-1,0,1\}$, there exist $c\in\F\setminus\{0\}$ and an orthogonal   matrix $O\in \mathcal{M}_n(\F)$ such that
\be{equation}{\label{Sn two map preserve traces power not -1}
\phi(A) = c^{-k} O AO^t,\quad \psi(B) = cOBO^t,
\quad  A, B\in \mathcal{S}_n(\F).
}
\end{enumerate}
\end{theorem}

\begin{proof} Assumption (2) leads to assumption (1) (cf. the proof of Theorem \ref{Two traces}). 
We prove the theorem under assumption (1). 
It suffices to prove the necessary part. %%Let   $\mathcal{P}_n(\R)$  denote the set of real   positive definite  matrices. 

Obviously, $\mathcal{S}_n\cap \mathcal{H}_n=\mathcal{S}_n(\R)$ and $\mathcal{S}_n\cap \mathcal{P}_n=\mathcal{P}_n(\R)$.
Let $S':=\{B\in \mathcal{P}_n(\R): B^{1/k}\in S\}$, which  is an open neighborhood of $I_n$ in $\mathcal{P}_n(\R)$ and whose real (resp. complex) span is $\mathcal{S}_n(\R)$ (resp. $\mathcal{S}_n$).
Using an analogous argument of the proof of Theorem \ref{Two traces}, and replacing $\mathcal{M}_n$ by $\mathcal{S}_n(\F)$,
replacing the basis \eqref{Mn projection basis} of $\mathcal{M}_n$   by
the following basis of  rank 1 projections  in $\mathcal{S}_n(\F)$:
\be{equation}{
\{E_{ii}: 1\le i\le n\}\cup \left\{\frac{1}{\sqrt{2}}(E_{ii}+E_{jj}+E_{ij}+E_{ji}):1\le i<j\le n\right\},
}
and replacing the usage of Theorem \ref{thm: Mn k-power negative}  by that of Theorem \ref{thm: Sn k-power negative},
we can prove the case $k\ge 2$, and for $k<0$, we can get the corresponding equalities up to \eqref{Mn T2 T1}.

Define a linear map $\psi_1:\mathcal{S}_n(\F)\to \mathcal{M}_n(\F)$  by $\psi_1(B):=\psi(B)\psi(I_n)^{-1}$.

When $k=-1$, we get the corresponding equality of \eqref{Mn psi B^2 identity}, so
that $\psi_1(B^2)=\psi_1(B)^2$ for $B\in \mathcal{S}_n(\F)$.
Similar to the proof of Theorem \ref{thm: Sn k-power negative}, we get $\psi_1(B^r)=\psi_1(B)^r$ for all $r\in \Z_+$.
By \cite[Theorem 6.5.3]{ZTC}, there is an invertible matrix $P\in \mathcal{M}_n(\F)$ such that
$\psi_1(B)=PBP^{-1}$, so that $\psi(B)=PBP^{-1}\psi(I_n)$ for $B\in \mathcal{S}_n(\F)$. Since $\psi(B)=\psi(B)^t$, we get
\be{equation}{
(P^{-1}\psi(I_n)P^{-t})B=B(P^{-1}\psi(I_n)P^{-t}),\quad  B\in \mathcal{S}_n(\F).
}
Therefore, $P^{-1}\psi(I_n)P^{-t}=cI_n$ for certain $c\in\F\setminus\{0\}$, so that
$\psi(B)=cPBP^t$ for all $B\in \mathcal{S}_n(\F)$. Consequently, we get \eqref{Sn two map preserve traces power -1}. The remaining claims are obvious.

When $k<-1$, using analogous argument as in the proof of $k<-1$ case of Theorem \ref{Two traces} and applying Theorem \ref{thm: Sn k-power negative},
we can get \eqref{Sn two map preserve traces power not -1}.
\end{proof}

%%%%%%%%%%%%%%%%%%%%%%%%%%%%%%%%%%%%%%%%%%%%%

\section{$k$-power linear preservers and trace of power-product preservers on $\mathcal{P}_n$ and $\mathcal{P}_n(\R)$}

In this section, we will determine $k$-power linear preservers and trace of power-product preservers on maps $\mathcal{P}_n\to \ol{\mathcal{P}_n}$
(resp. $\mathcal{P}_n(\R)\to \ol{\mathcal{P}_n(\R)}$).
Properties of such maps can be applied to   maps $\mathcal{P}_n\to \mathcal{P}_n$
and $\ol{\mathcal{P}_n}\to\ol{\mathcal{P}_n}$ (resp.  $\mathcal{P}_n(\R)\to \mathcal{P}_n(\R)$
and $\ol{\mathcal{P}_n(\R)}\to\ol{\mathcal{P}_n(\R)}$).

\subsection{$k$-power linear preservers on $\mathcal{P}_n$ and  $\mathcal{P}_n(\R)$}

\be{theorem}
{\label{thm: Pn linear map preserving k power}
Fix $k\in\Z\setminus\{0,1\}$.
A nonzero linear map $\psi:\mathcal{P}_n\to \ol{\mathcal{P}_n}$ (resp. $\psi:\mathcal{P}_n(\R)\to \ol{\mathcal{P}_n(\R)}$) satisfies that
\be{equation}{\label{Pn preserve k power}
\psi(A^k)=\psi(A)^k
}
on an open neighborhood of $I_n$ in $\mathcal{P}_n$ (resp. $\mathcal{P}_n(\R)$)
if and only if there is a unitary  (resp. real orthogonal) matrix $U\in M_n$ such that 
\be{equation}
{\label{Pn: trace power}
\psi(A)=U^*AU,\quad   A\in \mathcal{P}_n; \quad\text{or}\quad \psi(A)=U^*A^t U,\quad   A\in \mathcal{P}_n.
}
}

\begin{proof} We prove the case  $\psi:\mathcal{P}_n\to \ol{\mathcal{P}_n}$. The sufficient part is obvious. About the necessary part,
the nonzero linear map $\psi: \mathcal{P}_n\to \ol{\mathcal{P}_n}$  can be easily extended to a linear map $\wt{\psi}: \mathcal{H}_n\to \mathcal{H}_n$ that satisfies
$\wt{\psi}(A^k)=\wt{\psi}(A)^k$ on  an open neighborhood of $I_n$. By Theorem \ref{thm: Hn linear map preserving k power},
we immediately get \eqref{Pn: trace power}.

The case $\psi:\mathcal{P}_n(\R)\to \ol{\mathcal{P}_n(\R)}$ can be similarly proved using Theorem \ref{thm: Sn k-power negative}. 
\end{proof}

\subsection{Trace of powered product preservers on $\mathcal{P}_n$ and  $\mathcal{P}_n(\R)$}

Now consider the maps $\mathcal{P}_n\to\ol{\mathcal{P}_n}$ (resp. $\mathcal{P}_n(\R)\to\ol{\mathcal{P}_n(\R)}$) that
preserve trace of powered products.   
Unlike $\mathcal{M}_n$ and $\mathcal{H}_n$, the set  $ \mathcal{P}_n$ (resp. $ \mathcal{P}_n(\R)$) is not a vector space.
The trace of powered product preservers of two maps  have the following forms.

\begin{theorem}[Huang, Tsai \cite{Huang21}]\label{thm: Pn preserve trace power}
Let $a, b, c, d\in\R\setminus\{0\}$. Two maps $\phi,\psi:\mathcal{P}_n\to \ol{\mathcal{P}_n}$ satisfy
\begin{equation}\label{Pn trace power}
\tr (\phi(A)^{a} \psi(B)^{b})=\tr(A^c B^d),\quad   A, B\in \mathcal{P}_n,
\end{equation}
 if and only if   there exists an invertible $P\in \mathcal{M}_n$ such that
 \be{equation}{\label{Pn to Pn preserve trace power}
\be{cases}{
\phi(A)=(P^*A^c P)^{1/a} \\ \psi(B)=(P^{-1} B^d P^{-*})^{1/b}
}
  \text{ or } 
\be{cases}{
\phi(A)=[P^*(A^t)^{c}P]^{1/a}\\ \psi(B)=[P^{-1} (B^t)^{d} P^{-*}]^{1/b}
}
  A, B\in \mathcal{P}_n.
}
\end{theorem}

\begin{theorem}[Huang, Tsai \cite{Huang21}] \label{thm: multiple weighted product trace preserver}
Given an integer $m\ge 3$ and real numbers $\alpha_1,\ldots,\alpha_m,\beta_1,\ldots,\beta_m\in\R\setminus\{0\}$,
maps $\phi_i: \mathcal{P}_n\to \ol{\mathcal{P}_n}$ ($i=1,\ldots,m$) satisfy that
\begin{equation}\label{trace preserving multiple weighted products}
\tr \left(\phi_1(A_1)^{\alpha_1}\cdots \phi_m(A_m)^{\alpha_m}\right)=\tr \left(A_1^{\beta_1}\cdots A_m^{\beta_1}\right),\quad  A_1,\ldots,A_m\in \mathcal{P}_n,
\end{equation}
if and only if they have the following forms for certain $c_1,\ldots, c_m\in\R_+$ with $c_1\cdots c_m=1$:
\begin{enumerate}
\item When $m$ is odd, there exists a unitary matrix $U\in \mathcal{M}_n$ such that for $i=1,\ldots,m$:
\begin{equation}\label{trace preserving multiple weighted products odd}
\phi_i(A)=c_i^{1/\alpha_i} U^*A^{\beta_i/\alpha_i} U,\quad   A\in \mathcal{P}_n.
\end{equation}

\item When $m$ is even, there exists an invertible $M\in \mathcal{M}_n$ such that for $i=1,\ldots,m$:
\begin{equation}\label{trace preserving multiple weighted products even}
\phi_i(A)=
\begin{cases}
c_i^{1/\alpha_i} \left(M^*A^{\beta_i} M\right)^{1/\alpha_i}, &\text{$i$ is odd,}
\\
c_i^{1/\alpha_i} \left(M^{-1}A^{\beta_i} M^{-*}\right)^{1/\alpha_i} , &\text{$i$ is even,}
\end{cases}
\quad   A\in \mathcal{P}_n.
\end{equation}
\end{enumerate}
\end{theorem}

Both   Theorems \ref{thm: Pn preserve trace power} and \ref{thm: multiple weighted product trace preserver} can be analogously  extended to maps
$\mathcal{P}_n(\R)\to\ol{\mathcal{P}_n(\R)}$ without difficulties. 

Theorem \ref{thm: Pn preserve trace power} determines maps $\phi,\psi: \mathcal{P}_n\to\ol{\mathcal{P}_n}$ that satisfy \eqref{Pn trace power} throughout their domain.
If we  only assume the equality  \eqref{Pn trace power} for $(A,B)$ in  certain subset of $\mathcal{P}_n\times \mathcal{P}_n$ and assume certain linearity of $\phi$ and $\psi$, then $\phi$ and $\psi$ may have slightly different forms. We determine the case $a=c=1$ and $b=d=k\in\Z\setminus\{0\}$ here.

\begin{theorem}\label{thm: Pn preserve trace power open set}
 Let $k\in\Z\setminus\{0\}$.  Let $S$ be an open neighborhood of $I_n$ in $\mathcal{P}_n$.  Two maps $\phi,\psi:\mathcal{P}_n\to \ol{\mathcal{P}_n}$ satisfy
\begin{equation}\label{Pn trace power open set}
\tr (\phi(A)  \psi(B)^k )=\tr(A B^k),
\end{equation}
\begin{enumerate}
\item for all $A, B\in \mathcal{P}_n$, or
\item for all $A\in S,\ B\in \mathcal{P}_n$, and $\phi$ is linear,
\end{enumerate}
%% and that the map $\psi_1(B):=\psi(B^{1/d})^b$ is linear in $\mathcal{P}_n$,
if and only if there exists an invertible $P\in \mathcal{M}_n$ such that
  \be{equation}{\label{Pn to Pn preserve trace k power product}
\be{cases}{
\phi(A)=P^*A P \\ \psi(B)=(P^{-1} B^k P^{-*})^{1/k}
}
 \text{ or }
\be{cases}{
\phi(A)=P^*A^tP\\ \psi(B)=[P^{-1} (B^t)^{k} P^{-*}]^{1/k}
}
\quad  A, B\in \mathcal{P}_n.
}
The maps $\phi$ and $\psi$ satisfy \eqref{Pn trace power open set}
\begin{enumerate}
\setcounter{enumi}{2}
\item for all $A\in \mathcal{P}_n,\ B\in S$, and $\psi$ is linear, or
\item for all $A, B\in S$ and both $\phi$ and $\psi$ are  linear,
\end{enumerate}
if and only if when $k\in\{-1,1\}$, $\phi$ and $\psi$ take the form \eqref{Pn to Pn preserve trace k power product}, and
when $k\in\Z\setminus\{-1,0,1\}$, there exist a unitary matrix $U\in \mathcal{M}_n$ and $c\in\R^+$ such that
\be{equation}{\label{Pn two map preserve traces k not -1}
\be{cases}{
\phi(A) = c^{-k} U^*AU\\
\psi(B) = cU^*BU
}  \text{ or }
\be{cases}{
\phi(A) = c^{-k} U^*A^t U\\
\psi(B) = cU^*B^t U
}
  A, B\in \mathcal{P}_n.
}
\end{theorem}

\begin{proof} It suffices to prove the necessary part.

The case of assumption (1)  has been proved by  Theorem \ref{thm: Pn preserve trace power}.

Similar to  the proof of Theorem \ref{Two traces}, assumption (2) implies assumption (1); assumption (4)  implies assumption (3).
It remains to prove the case with assumption (3).

When $k=1$, assumption (3) is analogous to assumption (2), and we get \eqref{Pn to Pn preserve trace k power product}.

Suppose $k\in\Z\setminus\{1,0\}$.
Let $\psi_1:\mathcal{P}_n\to\ol{\mathcal{P}_n}$ be defined by $\psi_1(B):=\psi(B^{1/k})^{k}$. Let $\mathcal{S}_1:=\{B\in \mathcal{P}_n: B^{1/k}\in S\}.$ Then \eqref{Pn to Pn preserve trace k power product} with assumption (2) becomes
\be{equation}{
\tr(\phi(A)\psi_1(B))=\tr(AB),\quad  A\in \mathcal{P}_n,\ B\in \mathcal{S}_1.
}
Let $\wt{\psi}:\mathcal{H}_n\to \mathcal{H}_n$ be the linear extension of $\psi$.
By Theorem \ref{thm: two maps preserving trace},  $\phi$   can be extended to a linear bijection $\wt{\phi}:\mathcal{H}_n\to \mathcal{H}_n$ such that
\be{equation}{
\tr(\wt{\phi}(A)\wt{\psi}(B)^k)=\tr(\wt{\phi}(A)\psi_1(B^k))=\tr(AB^k),\quad  A\in \mathcal{H}_n,\ B\in S.
}
By Theorem \ref{thm: Hn Two traces} and taking into account the ranges of $\phi$ and $\psi$, we see that when $k=-1$, $\phi$ and $\psi$ take the form of
\eqref{Pn to Pn preserve trace k power product},
and
when $k\in\Z\setminus\{-1,0,1\}$,  $\phi$ and $\psi$ take the form of \eqref{Pn two map preserve traces k not -1}.
\end{proof}

Theorem \ref{thm: Pn preserve trace power open set} has counterpart results for $\phi,\psi:\mathcal{P}_n(\R)\to \ol{\mathcal{P}_n(\R)}$ and the proof is analogous using Theorem \ref{thm: Sn Two traces} instead of  Theorem \ref{thm: Hn Two traces}.

\section{$k$-power linear preservers and trace of power-product preservers on $\mathcal{D}_n$ and $\mathcal{D}_n(\R)$}

Let $\F=\C$ or $\R$.  Define the function
$\diag: \F^n\to \mathcal{D}_n(\F)$
to be the linear bijection  that sends each $(c_1,\cdots,c_n)^t$ to the diagonal matrix with $c_1,\ldots,c_n$ (in order) as the diagonal entries. Define
$\diag^{-1}:\mathcal{D}_n(\F)\to \F^n$
 the inverse map of $\diag$.

With the settings, every linear map $\psi: \mathcal{D}_n(\F)\to \mathcal{D}_n(\F)$ uniquely corresponds to a matrix $L_{\psi}\in \mathcal{M}_n(\F)$ such that
\be{equation}{\label{Dn linear map}
\psi(A)=\diag (L_{\psi} \diag^{-1}(A)),\quad  A\in \mathcal{D}_n(\F).
}

\subsection{$k$-power linear preservers on $\mathcal{D}_n$ and $\mathcal{D}_n(\R)$}

We define the linear functionals $f_i:\mathcal{D}_n(\F)\to\F$  $(i=0,1,\ldots,n)$, such that for each $A=\diag(a_1,\ldots,a_n)\in \mathcal{D}_n(\F)$,
\be{equation}{
f_0(A)=0;\quad  f_i(A)=a_i,\quad i=1,\ldots,n.
}

\begin{theorem}\label{thm: Dn k power}
Let $\F=\C$ or $\R$. Let $k\in\Z\setminus\{0,1\}$.  Let $S$ be an open neighborhood of $I_n$ in $\mathcal{D}_n(\F)$.  A   linear map $\psi: \mathcal{D}_n(\F)\to \mathcal{D}_n(\F)$ satisfies that
\be{equation}{\label{Dn k power preserving}
\psi(A^k)=\psi(A)^k,\quad  A\in S,
}
if and only if
\be{equation}{\label{Dn psi A form gen}
\psi(A)=\psi(I_n)\diag(  f_{p(1)}(A),\ldots, f_{p(n)}(A)),\quad  A\in \mathcal{D}_n(\F),
}
in which  $\psi(I_n)^k=\psi(I_n)$  and  $p:\{1,\ldots,n\}\to\{0,1,\ldots,n\}$ is a function such that
 $p(i)\ne 0$ when  $k<0$ for $i=1,\ldots,n$.
In particular, a linear bijection $\psi: \mathcal{D}_n(\F)\to \mathcal{D}_n(\F)$ satisfies \eqref{Dn k power preserving} if and only if
there is  a diagonal matrix $C\in \mathcal{M}_n(\F)$ with $C^{k-1}=I_n$ and a permutation matrix $P\in \mathcal{M}_n(\F)$ such that
\be{equation}{\label{Dn k power bijection}
\psi(A)= PCA P^{-1}, \quad  A\in \mathcal{D}_n(\F).
}
\end{theorem}

\begin{proof}
For every $A=\diag(a_1,\ldots,a_n)\in \mathcal{D}_n(\F)$, when $x\in\F$ is sufficiently close to $0$, we have $I_n+xA\in S$ and the power series of $(I_n+xA)^k$ converges, so that
$\psi((I_n+xA)^k)=\psi(I_n+xA)^k$.
\be{eqnarray}{
\psi((I_n+xA)^k) &=&  \psi(I_n)+xk\psi(A)+x^2\frac{k(k-1)}{2}\psi(A^2)+\cdots
\\
\psi(I_n+xA)^k &=&   \psi(I_n)+xk \psi(I_n)^{k-1}\psi(A)+x^2\frac{k(k-1)}{2}\psi(I_n)^{k-2}\psi(A)^2+\cdots
\quad
}
So for all $A\in \mathcal{D}_n(\F)$:
\be{eqnarray}{\label{Dn psi A identity}
\psi(A) &=& \psi(I_n)^{k-1}\psi(A),
\\\label{Dn psi A^2 identity}
\psi(A^2) &=& \psi(I_n)^{k-2}\psi(A)^2.
}
The linear map $\psi_1(A):= \psi(I_n)^{k-2}\psi(A)$ satisfies that
\be{equation}{\label{Dn psi_1 identity}
\psi_1(A^2)=\psi_1(A)^2,\quad  A\in \mathcal{D}_n(\F).
}
By \eqref{Dn linear map}, let $L_{\psi_1}=\mtx{\ell_{ij}}\in \mathcal{M}_n(\F)$ such that
$\diag^{-1}(\psi_1(A))=L_{\psi_1}(\diag^{-1}(A))$ for $A\in \mathcal{D}_n(\F)$. Then
\eqref{Dn psi_1 identity} implies that for all $A=\diag(a_1,\ldots,a_n)\in \mathcal{D}_n(\F)$:
\be{equation}{
\sum_{j=1}^{n} \ell_{ij}a_j^2 = \left(\sum_{j=1}^{n} \ell_{ij}a_j\right)^2,\quad  i=1,2,\ldots,n.
}
Therefore, each row of $L_{\psi_1}$ has at most one nonzero entry and each nonzero entry must be  1.
We get
\be{equation}{\label{Dn psi_1}
\psi_1(A)=\diag(f_{p(1)}(A),\ldots,f_{p(n)}(A))
}
in which $p:\{1,\ldots,n\}\to\{0,1,\ldots,n\}$ is a function. Suppose $\psi(I_n)=\diag(\lambda_1,\ldots,\lambda_n)$. Then
\eqref{Dn psi A identity} implies that $\psi(A)=\psi(I_n)\psi_1(A)$ has the form \eqref{Dn psi A form gen}. Obviously, 
$\psi(I_n)^k=\psi(I_n)$ and
when $k<0$,   each $p(i)\ne 0$ for $i=1,\ldots,n$. Moreover,
when $\psi$ is a linear bijection, \eqref{Dn psi_1} shows that $\psi_1(A)=PAP^{-1}$ for a permutation matrix $P$. 
\eqref{Dn k power bijection} can be easily derived. 
\end{proof}

\subsection{Trace of power-product preservers on $\mathcal{D}_n$ and $\mathcal{D}_n(\R)$}

In \cite{Huang21}, we show that two maps $\phi,\psi:\mathcal{D}_n(\F)\to \mathcal{D}_n(\F)$ satisfy
$\tr(\phi(A)\psi(B))=\tr(AB)$ for $A,B\in \mathcal{D}_n(\F)$ if and only if  there exists an invertible $N\in \mathcal{M}_n(\F)$ such that
\be{equation}
{
\phi(A)=\diag(N\diag^{-1}(A)),\ \psi(B)=\diag(N^{-t}\diag^{-1}(B)),\quad  A, B\in \mathcal{D}_n(\F).
}
When $m\ge 3$, the maps $\phi_1,\ldots,\phi_m: \mathcal{D}_n(\F)\to \mathcal{D}_n(\F)$ satisfying
$\tr(\phi_1(A_1)\cdots\phi_m(A_m))=\tr(A_1\cdots A_m)$ for $A_1,\ldots,A_m\in \mathcal{D}_n(\F)$ are also determined in \cite{Huang21}.

Next we consider the trace of power-product preserver on $\mathcal{D}_n(\F)$.

\begin{theorem}\label{thm: Dn Two traces2} 
Let $\F=\C$ or $\R$.  Let $k\in\Z\setminus\{0,1\}$. Let $S$ be an open neighborhood of $I_n$ in $\mathcal{D}_n(\F)$.  Two maps $\phi,\psi: \mathcal{D}_n(\F)\to \mathcal{D}_n(\F)$   satisfy that
\be{equation}
{\label{Dn: two trace preserving}
\tr (\phi(A)\psi(B)^k)=\tr (AB^k), %%\quad  A \in \mathcal{D}_n(\F), \ B\in S,
}
\begin{enumerate}
\item for all $A\in \mathcal{D}_n(\F),\ B\in S$, and $\psi$ is linear, or
\item for all $A, B\in S$ and both $\phi$ and $\psi$ are  linear,
\end{enumerate}
if and only if there exist  an invertible diagonal matrix $C\in \mathcal{D}_n(\F)$ and a permutation matrix $P\in \mathcal{M}_n(F)$ such that
\be{equation}
{\label{Dn two map preserve traces power}
\phi(A)=PC^{-k}A P^{-1},\quad \psi(B)= PCB P^{-1}, \quad   A,B\in \mathcal{D}_n(\F).
}
\end{theorem}

\begin{proof} Assumption (2) leads to assumption (1) (cf. the proof of Theorem \ref{Two traces}). 
We prove the theorem under assumption (1). 

For every $B\in \mathcal{D}_n(\F)$,  $I_n+xB\in S$ and the power series of $(I_n+xB)^k$ converges when $x\in\F$ is sufficiently close to $0$, so that
\be{equation}{
\tr (\phi(A)\psi(I_n+xB)^k)=\tr (A(I_n+xB)^k)
}
Comparing degree one terms and degree two terms in the power series of the above equality, respectively, we get the following equalities for $A, B\in \mathcal{D}_n(\F)$: 
\be{eqnarray}{\label{Dn phi A psi B}
\tr(\phi(A)\psi(B)\psi(I_n)^{k-1}) &=& \tr(AB),
\\\label{Dn phi A psi B^2}
\tr(\phi(A)\psi(B)^2\psi(I_n)^{k-2}) &=&\tr(AB^2).
}
Applying Theorem \ref{thm: two maps preserving trace} to \eqref{Dn phi A psi B},  
$\psi(I_n)$ is invertible and
both $\phi$ and $\psi$ are linear bijections. \eqref{Dn phi A psi B} and \eqref{Dn phi A psi B^2} 
imply that $ \psi(B^2)\psi(I_n)^{k-1}=\psi(B)^2\psi(I_n)^{k-2}$.
Let  $\psi_1(B):=\psi(B)\psi(I_n)^{-1}$. Then $\psi_1(B^2)=\psi_1(B)^2$ for $B\in \mathcal{D}_n(\F)$.  
By   Theorem \ref{thm: Dn k power} and  $\psi_1(I_n)=I_n$,
 there exists a permutation matrix $P\in \mathcal{M}_n(F)$ such that
$\psi_1(B)=PBP^{-1}$ for $B\in \mathcal{D}_n(\F)$. So $\psi(B)= \psi(I_n)PBP^{-1}=PCBP^{-1}$ for $C:=P^{-1}\psi(I_n)P\in \mathcal{D}_n(\F)$. 
Then \eqref{Dn: two trace preserving} implies \eqref{Dn two map preserve traces power}.
\end{proof}

\section{$k$-power injective linear preservers and trace of power-product preservers on $\mathcal{T}_n$ and $\mathcal{T}_n(\R)$}

\subsection{$k$-power preservers on $\mathcal{T}_n(\F)$}

The characterization of injective linear $k$-power preserver on $\mathcal{T}_n(\F)$ can be derived from Cao and Zhang's   characterization of injective additive $k$-power preserver on $\mathcal{T}_n(\F)$ (\cite{Cao05} or \cite[Theorem 6.5.2]{ZTC}).

\begin{theorem} [Cao and Zhang \cite{Cao05}] \label{thm: Tn injective k power preserver}
 Let $k\geq2$ and $n\geq3$. Let $\F$ be a field with $\operatorname{char}(\F)=0$ or $\operatorname{char}(\F)>k$.  Then $\psi: \mathcal{T}_n(\F)\mapsto \mathcal{T}_n(\F)$  is an injective linear map such that $\psi(A^k)=\psi(A)^k$ for all $A\in \mathcal{T}_n(\F)$ if and only if there exists a $(k-1)$th root of unity $\lambda$ and an invertible matrix $P\in \mathcal{T}_n(\F)$ such that   
\be{eqnarray} { \label{linear bijection preserving power 4}
\psi(A)=\lambda PAP^{-1}, && A\in \mathcal{T}_n(\F),\quad\text{or}
\\ \label{linear bijection preserving power 5}
\psi(A)=\lambda PA^{-} P^{-1}, && A\in \mathcal{T}_n(\F),
 }
 where $A^{-}=(a_{n+1-j,n+1-i})$ if $A=(a_{ij})$.
\end{theorem}

\begin{example}
When $n=2$, the injective linear maps that  satisfy
$\psi(A^k)=\psi(A)^k$ for $A\in  \mathcal{T}_2(\F)$ send $A=\mtx{a_{11} &a_{12}\\0 &a_{22}}$ to  the following $\psi(A)$:
\be{equation}{
\lambda\mtx{a_{11} &c a_{12}\\0 &a_{22}},\quad
\lambda\mtx{a_{22} &c a_{12}\\0 &a_{11}},
}
in which $\lambda^{k-1}=1$ and $c\in\F\setminus\{0\}.$ 
\end{example}

\begin{example} Theorem \ref{thm: Tn injective k power preserver} does not hold if $\psi$ is not assumed to be injective. 
Let $n=3$ and suppose  $\psi:\mathcal{T}_3(\F)\to \mathcal{T}_3(\F)$ is a linear map that  sends $A=\mtx{a_{ij}}_{3\times 3}\in\mathcal{T}_3(\F)$ to one of the following $\psi(A)$ ($c,d\in\F$): 
\be{equation}{
\mtx{a_{11} &c a_{12} &0\\0 &a_{22} &d a_{23}\\0 &0 &a_{33}}, 
%\mtx{a_{11} &a_{12} &0\\0 &a_{22} &0 \\0 &0 &a_{33}},
\mtx{a_{33} &0 &0\\0 &a_{11} &0 \\0 &0 &a_{22}}, 
\mtx{a_{22} &0 &c a_{12}\\0 &0 &0 \\0 &0 &a_{11}}.
}
Then each $\psi$ satisfies that $\psi(A^k)=\psi(A)^k$ for every positive integer $k$ but it is not of the forms in Theorem \ref{thm: Tn injective k power preserver}.
\end{example}

We extend Theorem \ref{thm: Tn injective k power preserver} to the following result that includes negative $k$-powers and that
only assumes $k$-power preserving in a neighborhood of $I_n$. 

\begin{theorem} \label{thm: Tn injective k-power negative}
Let $\F=\C$ or $\R$. Let  integers $k\ne 0, 1$ and $n\ge 3$.
Suppose that $\psi: \mathcal{T}_n(\F)\to \mathcal{T}_n(\F)$ is an injective linear map such that $\psi(A^k)=\psi(A)^k$ for all $A$ in an open neighborhood of $I_n$ in $\mathcal{T}_n(\F)$ consisting of invertible matrices. Then there exist   $\lambda\in\F$ with $\lambda^{k-1}=1$ and an invertible matrix $P\in \mathcal{T}_n(\F)$ such that
\be{eqnarray}
{\label{Tn preserving k-power 1}
\psi(A)=\lambda PAP^{-1}, &&  A\in \mathcal{T}_n(\F),\quad\text{or}
\\ \label{Tn preserving k-power 2}
\psi(A)=\lambda PA^{-}P^{-1}, &&   A\in \mathcal{T}_n(\F).
}
 where $A^{-}=(a_{n+1-j,n+1-i})=J_n A^t J_n$ if $A=(a_{ij})$, $J_n$ is the anti-diagonal identity.
\end{theorem}

\begin{proof}
Obviously $\psi$ is a linear bijection. Follow  the same process in the proof of Theorem 
\ref{thm: Mn k-power negative}. In both $k\ge 2$ and $k<0$ cases we have
$\psi(I_n)$ commutes with the range of $\psi$, so that $\psi(I_n)=\lambda I_n$ for 
$\lambda\in\F$ and $\lambda^{k-1}=1.$  Moreover, let $\psi_1(A):=\psi(I_n)^{-1}\psi(A)$,
then $\psi_1$ is injective linear and $\psi_1(A^2)=\psi_1(A)^2$ for $A\in\mathcal{T}_n(\F).$ 
Theorem 
\ref{thm: Tn injective k power preserver} shows that $\psi_1(A)=PAP^{-1}$ or $\psi_1(A)=PA^{-}P^{-1}$ 
for certain invertible $P\in \mathcal{T}_n(\F).$
It leads to \eqref{Tn preserving k-power 1} and \eqref{Tn preserving k-power 2}. 
\end{proof}

\subsection{Trace of power-product preservers on $\mathcal{T}_n$ and $\mathcal{T}_n(\R)$}

Theorem \ref{thm: two maps preserving trace} or 
 Corollary \ref{thm: trace preserver on special spaces} 
does not work for maps on $\mathcal{T}_n(\F)$. 
However, the following trace preserving result can be easily derived from  Theorem \ref{thm: Dn Two traces2}. We have $\mathcal{T}_n(\F)=\mathcal{D}_n(\F)\oplus \mathcal{N}_n(\F)$.   
 Let $\operatorname{D}(A)$ denote the diagonal matrix that takes the diagonal of $A\in \mathcal{T}_n(\F)$.

\begin{theorem}\label{Tn Two traces}
Let $\F=\C$ or $\R$.
Let  $k \in\Z\setminus\{0,1\}$. Let $S$ be an open neighborhood of $I_n$ in $\mathcal{T}_n(\F)$  consisting of invertible matrices.
Then two maps $\phi, \psi: \mathcal{T}_n(\F)\to \mathcal{T}_n(\F)$ satisfy that
\be{equation}
{\label{Tn trace two maps power}
 \tr (\phi(A)\psi(B)^k)=\tr (AB^k),%\quad  A \in \mathcal{M}_n(\F),   B\in S,
}
\begin{enumerate}
\item for all $A\in \mathcal{T}_n(\F),\ B\in S$, and $\psi$ is linear, or
\item for all $A, B\in S$ and both $\phi$ and $\psi$ are  linear,
\end{enumerate}
if and only if  $\phi$ and $\psi$     send  $\mathcal{N}_n(\F)$ 
to $\mathcal{N}_n(\F)$, $(\operatorname{D}\circ\phi)|_{\mathcal{D}_n(\F)}$  and $(\operatorname{D}\circ\psi)|_{\mathcal{D}_n(\F)}$
are linear bijections characterized by \eqref{Dn two map preserve traces power} in Theorem \ref{thm: Dn Two traces2}, and
$\operatorname{D}\circ\phi=\operatorname{D}\circ\phi\circ \operatorname{D}$. 
\end{theorem}

\begin{proof} The sufficient part is easy to verify. We prove the necessary part here. 
Let $\phi':=(\operatorname{D}\circ\phi)|_{\mathcal{D}_n(\F)}$ and $\psi':=(\operatorname{D}\circ\psi)|_{\mathcal{D}_n(\F)}$. Then $\phi',\psi': \mathcal{D}_n(\F)\to \mathcal{D}_n(\F)$ satisfy 
$\tr(\phi'(A)\psi'(B)^k)=\tr(AB^k)$ for $A,B\in \mathcal{D}_n(\F)$. So they are characterized by \eqref{Dn two map preserve traces power}. The bijectivity of $\phi'$ and $\psi'$ implies that 
 $\phi$ and $\psi$ must send $\mathcal{N}_n(\F)$ 
to $\mathcal{N}_n(\F)$ in order to satisfy \eqref{Tn trace two maps power}. Moreover,
$\phi$ should send matrices with same diagonal to matrices with same diagonal, which implies that $\operatorname{D}\circ\phi=\operatorname{D}\circ\phi\circ \operatorname{D}$. 
\end{proof}

\end{document}